\documentclass[12pt,reqno]{amsart}

\usepackage[utf8]{inputenc}
\usepackage[english]{babel}
\usepackage{newlfont}
\usepackage{color}
\usepackage{ amssymb }
\usepackage{graphicx}
\usepackage{amsmath}
\usepackage{scalerel}
\usepackage{stackengine,wasysym}
\usepackage{mathrsfs}
\usepackage{enumitem}
\usepackage{braket}
\usepackage{blindtext}
\usepackage[utf8]{inputenc}
\usepackage{fancyhdr,appendix,psfrag,layout}
\usepackage{hyperref}
\usepackage{url}
\usepackage[capitalise,nameinlink,noabbrev]{cleveref}
\crefname{equation}{}{}
\crefname{enumi}{}{} 
\usepackage{amsfonts}
\usepackage{amsthm}
\usepackage[english]{babel}
\usepackage{fancyhdr}
\usepackage{tikz}
\usepackage{array}
\usepackage{multirow}
\usetikzlibrary{positioning}
\addto\captionsenglish{}

\newtheorem{lemma}{Lemma}[section]

\newtheorem{theorem}[lemma]{Theorem}

\theoremstyle{definition}

\newtheorem{remark}{Remark}

\numberwithin{equation}{section}

\def\OO{{\mathcal O}}
\def\BB{{\mathcal B}}

\def\O{{\mathrm{O}}}
\def\SO{\mathrm{SO}}

\def\PSO{\mathrm{PSO}}

\def\In{\mbox{Invol}}
\def\I{\mbox{I}}
\def\Rac{\mbox{Rac}}
\def\mod{{\mathrm{mod}}}
\def\Dih{\mbox{Dih}}

\def\Isompr{\mbox{Isom}\!^+\!}

\def\Norm{\mbox{Norm}}

\def\C{{\mathbb C}}
\def\R{{\mathbb R}}
\def\H{{\mathbb H}}
\def\Z{{\mathbb Z}}

\begin{document}

\textwidth=450pt\oddsidemargin=0pt

\title{Involutions of spherical 3-manifolds}

\author[M. Mecchia]{Mattia Mecchia}
\address{M. Mecchia: Dipartimento Di Matematica e Geoscienze, Universit\`{a} degli Studi di Trieste, Via Valerio 12/1, 34127, Trieste, Italy.} \email{mmecchia@units.it}

\author[B. Schilling]{Baptiste Schilling}
\address{B. Schilling: Institut Camille Jordan, Universit\'e Claude Bernard Lyon 1, Bâtiment Jean Braconnier, 21 avenue Claude Bernard, 69100 Villeurbanne, France} \email{schilling@math.univ-lyon1.fr}




\begin{abstract}
We classify involutions acting on spherical 3-manifolds up to conjugacy. Our geometric approach provides insights into numerous topological properties of these involutions.

\end{abstract}

\maketitle
  
\section{Introduction}

We call \textit{involution} a periodic diffeomorphism of order 2 acting on a smooth manifold. This paper considers orientation-preserving involutions acting on orientable and closed spherical 3-manifolds. A closed and orientable spherical 3-manifold is the quotient of $S^3$ by a finite subgroup of $\SO(4)$ acting freely on $S^3$.  Two spherical 3-manifolds $S^3/G$ and $S^3/G'$ are homeomorphic if and only if $G$ and $G'$ are conjugate in $\O(4).$ Thus, the classification of the spherical 3-manifolds is equivalent to the classification up to conjugacy of the finite subgroups acting freely on $S^3$, which is originally due to Seifert and Threlfall (\cite{threlfall-seifert} and \cite{threlfall-seifert2}.)  If the finite group $G$ does not act freely the quotient $S^3/G$ has a natural structure as smooth spherical 3-orbifold. Also in this case, the quotient orbifolds are diffeomorphic if and only if the corresponding groups are conjugate. The classification up to conjugacy of the finite subgroups of $\SO(4)$ is summarised in Table~\ref{subgroup}; the groups are divided into several families. Section~\ref{sec-subgroups} explains which subgroups act freely on $S^3$; in particular,  free actions can occur only for Families 1, $1'$, 2, 3, 5, 6, 7 and 9.

The study of involutions acting on 3-manifolds is a classical research theme in the field of 3-manifolds. This topic played a significant role in the development of Thurston's geometrization program. For spherical 3-manifolds, the cases of lens spaces and quaternion manifolds have been considered in several papers, see for example \cite{Hodgson-Rubinstein, myers, kim, Rubinstein}.   Following the proof of Thurston's geometrization theorem, the analysis of the involutions of geometric 3-manifolds is easier since we know that each involution is conjugate to an isometry and the whole analysis can be carried out in the isometry groups of the manifolds.  However, a comprehensive description of the involutions of spherical 3-manifolds is absent from the literature. This paper aims to fill that gap by providing a detailed classification of involutions for all spherical 3-manifolds, along with geometric insights into their actions. Specifically, we obtain the algebraic classification of the involutions up to conjugacy and describe the corresponding quotient orbifolds. In particular, we can distinguish if an involution acts freely,  or equivalently if the quotient orbifold is a manifold. If the quotient orbifold has $S^3$ as underlying topological space the involution is called \textit{hyperelliptic}. The hyperelliptic involutions have been largely studied in dimensions two and three. In the 3-dimensional case if a manifold admits a hyperelliptic involution it can be seen as the 2-fold branched cover of $S^3$ along a link that is the singular set of the quotient orbifold. In our case, we can individuate which involutions are hyperelliptic. 

The analysis of the quotients can be carried out by using Seifert fibrations.   Seifert introduced this structure for 3-manifolds, see \cite{seifert}. Every spherical 3-manifold admits such a fibration. Bonahon and Siebenmann \cite{bonahon-sibenmann} generalized the definition to 3-orbifolds. There exist spherical 3-orbifolds not admitting a Seifert fibration, but any orbifold obtained taking the quotient of a spherical manifold by an involution is Seifert fibered.  A Seifert fibration is uniquely determined by a series of invariants (base orbifolds, local invariants, Euler class and boundary invariants) and, in our case, they can be computed by using the formulae given in \cite{mecchia-seppi}. For notations and details about Seifert fibrations, we refer to \cite{MecchiaSeppi3}.

The classification of involutions and their principal geometric features are summarised in Table~\ref{general quotient} and \ref{small index quotient}. The spherical 3-manifold $S^3/G$ is represented through the group $G$ in the first column. The conjugacy classes of involution of $S^3/G$ are represented by extensions of index 2 of $G$ that are listed in the second column. In the third column, it is said if the action of the involutions is free. In Table~\ref{general quotient}, if the extensions are in Family 1 or $1'$ we cannot produce a closed formula to describe when the involution acts freely. However, it is possible to use Tables~\ref{seifert-fam-1} and \ref{seifert-fam-1-prime} to analyse this aspect. In the fourth column, in the fourth column we say if the involution is hyperelliptic.

The spherical 3-manifolds are divided into some subclasses: lens spaces, prism manifolds and platonic manifolds, which are further subdivided into tetrahedral, octahedral and icosahedral manifolds.

If the group $G$ is in Families 1 and $1'$, the manifold  $S^3/G$ is a lens space; the groups in Families 2 and 3 give prism manifolds; tetrahedral manifolds are the quotients by groups in Families 5 and 6; finally, if $G$ is in Family 7 (resp. Family 8) the manifold $S^3/G$ is octahedral (resp. icosahedral).

The following statement concisely describes the classification according to the manifold type. 

\begin{theorem}\label{recap}
Let $M$ be a spherical 3-manifold.
\begin{enumerate}
    \item If $M$ is a lens space, it has three to eight involution classes. 
    \item  If $M$ is a prism manifold, it has three to five involution classes. Moreover, $M$ admits an involution acting freely if and only if it is diffeomorphic to $S^3/G$ where $G$ is in Family 2; such involution is unique up to conjugacy. 
    \item If $M$ is a tetrahedral manifold, it has two to five involution classes.  Moreover, $M$ admits an involution acting freely if and only if it is diffeomorphic to $S^3/G$ where $G$ is in Family 5; such involution is unique up to conjugacy.
   \item  If $M$ is an octahedral or an icosahedral manifold, it has one or two involution classes; no involution acts freely.
\end{enumerate}

Moreover, $M$ admits, up to conjugacy, exactly one hyperelliptic involution.

\end{theorem}

Lens spaces involutions were considered in \cite{myers, Hodgson-Rubinstein}; in particular, an explicit classification of involutions not acting freely is given in \cite{Hodgson-Rubinstein}.  Considering lens spaces, it is also worth mentioning the paper of Kalliongis and Miller \cite{kalliongis-miller} where the isometry groups of the lens spaces are described in many details. The involutions of prism manifolds are studied in \cite{Rubinstein} and \cite{kim}. 

The existence of a hyperelliptic involution for each manifold and its uniqueness up to conjugacy implies that each spherical 3-manifold can be represented in exactly one way as the 2-fold branched cover of $S^3$ along a link. We note that not all the manifolds can be obtained as the 2-fold branched cover of a link, and, if it exists,  this representation is generally highly not unique (see, for example, \cite{viro, mecchia-zimmermann, Kawauchi}). In Subsection~\ref{quotient}, we prove that in the spherical case, all the singular sets of quotient orbifolds obtained by hyperelliptic involutions are Montesinos links, so we can obtain the following result. 

\begin{theorem}\label{montesinos-link}
Each spherical 3-manifold is the 2-fold branched cover of $S^3$ along a Montesinos link with at most three tangles. If two links give diffeomorphic spherical 2-fold branched covers of $S^3$, they are equivalent.  
\end{theorem}

Montesinos links are introduced in \cite{montesinos}, see also \cite{burde-zieschang}; a brief description of this class of links is given in Subsection~\ref{quotient} where the proof of Theorem~\ref{montesinos-link} is completed. 

 Theorem~\ref{montesinos-link} can also be proved on the basis of the results already present in the literature and not using our classification:  Dunbar \cite{dunbar} gives a list of spherical orbifolds whose underlying topological space is  $S^3$.  We restrict our attention to the orbifolds whose singular set is a link. Many of these links are already presented as Montesinos links. There are  few exceptions but also these can be manipulated to obtain  Montesinos links.  The 2-fold branched cover of $S^3$ along a Montesinos link is described in \cite{montesinos} via a Seifert fibration.  A manifold might admit inequivalent Seifert fibrations but the topological classification of Seifert  3-manifolds up to homeomorphism is easy to handle except for lens spaces that admit infinite inequivalent fibrations, see \cite{orlik}. For lens spaces, the uniqueness is obtained in \cite{Hodgson-Rubinstein}. In the other cases, we can use the topological classification of  Seifert 3-manifolds to obtain the result. In the present paper, Theorem~\ref{montesinos-link} is a by-product of the classification.

Even if we skip the hypothesis of being hyperelliptic, our approach can obtain more geometric information; in Subsection~\ref{quotient}, some cases are analyzed. Finally, we note that our method can also be used to study the involutions of the spherical 3-orbifolds. A general classification would involve a huge number of cases, so it may not be worth the effort, but our approach can be easily extended to analyze single cases also in the orbifold setting.






\section{Finite subgroup of SO(4)}\label{sec-subgroups}

 We use the quaternion algebra to classify the finite subgroups of $\SO(4)$. We essentially follow \cite{duval} although it must be mentioned that in Du Val's list of finite subgroups of $\SO(4)$ there are three missing cases, see also \cite{conway-smith, mecchia-seppi, mecchia-seppi2}.

Let  us identify $\R^4$ with the quaternion algebra $\H=\{a+bi+cj+dk\,|\,a,b,c,d\in\R\}=\{z_1+z_2j\,|\,z_1,z_2\in\C\}$. Given $q=z_1+z_2 j\in\H$, its conjugate is $\bar q=\bar z_1-z_2 j$. Thus the standard positive definite quadratic form of $\R^4$ is identified as $q\bar q=|z_1|^2+|z_2|^2$.
The three-sphere $S^3$ is represented as  the set of unit quaternions:
\begin{equation*}
S^3=\{a+bi+cj+dk \,|\, a^2+b^2+c^2+d^2=1\}=\{z_1+z_2j\, \,|\, |z_1|^2+|z_2|^2=1\}\,.
\end{equation*}
which is thus endowed with a multiplicative group structure induced from that of $\H$.

The finite subgroups of $S^3$ are well known. Up to conjugacy, those are: 
\begin{itemize}
\item $C_n=\{e^{2i\pi k/n}, 0\leq k<n\}$
\item $D^*_{4n}= C_{2n}\cup C_{2n} j$
\item $T^*= D_8^*\cup \alpha D_8^* \cup \alpha^2 D_8^*$, where $\alpha= \dfrac{1}{2}(1+i+j+k)$
\item $O^* = T^*\cup \beta T^*$, where $\beta =\dfrac{1}{\sqrt{2}}(1+j)$
\item $I^* = \bigcup_{r=0}^4 \gamma^r T^*$, where $\gamma= \dfrac{1}{2}(\dfrac{1+\sqrt{5}}{2}^{-1}+ \dfrac{1+\sqrt{5}}{2}j +k)$
\end{itemize}

The group $C_n$ is a cyclic group of order $n$. In the following, we will write $C_n$ to refer to the subgroup of $S^3$, and $\mathbb{Z}_n$ to designate the abstract cyclic group of order $n$.

The group $D^*_{4n}$ is the binary dihedral group of order $4n$: It is a central extension of the dihedral group $D_{2n}$ by a group of order $2$, and admits the presentation $\langle \zeta, j | \zeta^{2n}=j\zeta j^{-1}\zeta= \zeta^n j^2=1 \rangle$. Observe in particular that $D_4^*=\{\pm 1, \pm j\}$ is conjugate to $C_4$ (for this reason, we will only refer to $D^*_{4n}$ when $n\geq 2$), and $D^*_8= \{\pm 1, \pm i, \pm j, \pm k \}$ is also called the quaternion group.

The groups $T^*, O^*$ and $I^*$ are respectively called binary tetrahedral, binary octahedral, and binary icosahedral groups, and their orders are $24, 48$, and $120$.

Let $\Phi: S^3\times S^3 \rightarrow SO(4)$ be the application given by $\Phi(p,q)(x)=pxq^{-1}.$ It is a $2:1$ surjective morphism whose kernel is equal to $\{\pm (1,1)\}$.
If $G$ is a subgroup of $SO(4)$, we write $\tilde{G}=\Phi^{-1}(G)$. Since $(-1,-1)$ is central, two subgroups $G$ and $G'$ are conjugate in $SO(4)$ if and only if $\tilde{G}$ and $\tilde{G'}$ are conjugate in $S^3\times S^3$, thus reducing the problem of classifying the finite subgroups of $SO(4)$ up to conjugacy to the classification of the finite subgroups of $S^3\times S^3$ containing $(-1,-1)$ up to conjugacy. 

Let $\tilde{G}$ be a finite subgroup of $S^3\times S^3$, and $pr_1, pr_2$ the projections on the first and second coordinate, respectively. We introduce the groups $L=pr_1(\tilde{G}),$ $L_K=pr_1(\tilde{G}\cap (S^3 \times \{1\}) ),$ $R=pr_2(\tilde{G})$ and $R_K=pr_2(\tilde{G}\cap (\{1\} \times S^3) )$.  The projections induce two isomorphisms $\overline{pr_1}: \tilde{G}/L_K\times R_K \rightarrow L/L_K$ and $\overline{pr_2}: \tilde{G}/L_K\times R_K \rightarrow R/R_K$, and we define an isomorphism  $\phi: L/L_K \rightarrow R/R_K$  by $\phi= \overline{pr_2}\circ \overline{pr_1}^{-1}$ .

Conversely, if $L$ and $R$ are finite subgroups of $S^3$ with normal subgroups $L_K$ and $R_K$ and an isomorphism $\phi: L/L_K \rightarrow R/R_K$, then there exists a unique subgroup $H$ of $S^3\times S^3$ such that $L=pr_1(H),$ $L_K=pr_1(H\cap (S^3 \times \{1\}) ),$ $R=pr_2(H),$ $R_K=pr_2(H\cap (\{1\} \times S^3) )$ and the  isomorphism $\phi$ is equal to $\overline{pr_2}\circ \overline{pr_1}^{-1}$. The definition of $H$ is the following: 
$$H=\{ (x,y)\in L\times R \,|\, \phi (xL_K)= yR_K \}.$$  We use the quintuple $(L,L_K, R,R_K, \phi)$ to denote the group $H.$

This presentation describes the subgroups of $S^3\times S^3$ in terms of subgroups of $S^3$. With the help of the following lemma, whose proof is straightforward, we can also understand if two subgroups of $S^3\times S^3$ are conjugate.

\begin{lemma}\label{conjugation}Let $\tilde{G}=(L, L_K , R, R_K , \phi)$ and $\tilde{G'} = (L' , L_K' , R'  , R_K' , \phi' )$ be finite subgroups
of $S^3\times S^3$ containing $(-1,-1)$. An element $(g, f )\in S^3\times S^3$ conjugates $\tilde{G}$ to $\tilde{G'}$ if and only if the following conditions are satisfied:
\begin{enumerate}
\item $g^{-1}Lg = L'$ and $f^{-1} R f = R'$
\item $g^{-1}L_Kg = L_K'$ and $f^{-1} R_K f = R_K'$
\item The equality $\bar{c}_f \circ \phi = \phi' \circ \bar{c}_g$ holds, i.e. the following diagram commutes:

\begin{center}

\begin{tikzpicture}
\node (A) {$L/L_K$};
\node (B) [below= of A]  {$L'/L'_K$};
\node (C) [right= of A] {$R/R_K$};
\node (D) [below= of C] {$R'/R_K'$};
\node (E) [left= of A] {$L$};
\node (F) at (-1.96, -1.71) {$L'$};
\node (G) [right= of C] {$R$};
\node (H) at (4.39, -1.71) {$R'$};

\draw[->] (A) -- (B) node [midway, left] {$ \bar{c}_g$};
\draw[->] (C) -- (D) node [midway, left] {$ \bar{c}_f$};
\draw[->] (A) -- (C) node [midway, above] {$ \phi$};
\draw[->] (B) -- (D) node [midway, above] {$ \phi'$};
\draw[->] (E) -- (A) node [midway, left] {};
\draw[->] (F) -- (B) node [midway, left] {};
\draw[->] (G) -- (C) node [midway, left] {};
\draw[->] (H) -- (D) node [midway, left] {};
\draw[->] (E) -- (F) node [midway, left] {$ c_g$};
\draw[->] (G) -- (H) node [midway, left] {$ c_f$};

\end{tikzpicture}

\end{center}
\end{enumerate}
\noindent
Where $\bar{c}_g$ (resp. $\bar{c}_f$) is the application $L/L_K \rightarrow L'/L_K'$ (resp. $R/R_K \rightarrow R'/R_K'$) induced in the quotient by $c_g: L\rightarrow L'$ (resp. $c_f: R\rightarrow R'$), the conjugation by $g$ (resp. by $f$).

\end{lemma}

Using this lemma and the list of subgroups of $S^3$, it is possible to obtain a complete classification up to conjugacy of the finite subgroups of $S^3\times S^3$ containing $\pm (1,1)$, reported in Table \ref{subgroup}.

\begin{table}\caption{Finite subgroups of $SO(4)$}\label{subgroup}
\begin{tabular}{|l|c|c|c|}
\hline
 & family of groups & order of $\Phi(G)$ & \\
\hline
 1. & $(C_{2mr}/C_{2m},C_{2nr}/C_{2n})_s$ & $2mnr$ & $\gcd(s,r)=1$ \\  
 $1^{\prime}$. & $(C_{mr}/C_{m},C_{nr}/C_{n})_s$ & $(mnr)/2$ & $\gcd(s,r)=1$ $\gcd(2,n)=1$ \\
&&& $\gcd(2,m)=1$  $\gcd(2,r)=2$\\
 2. & $(C_{2m}/C_{2m},D^*_{4n}/D^*_{4n})$ & $4mn$ &  \\ 
 3. & $(C_{4m}/C_{2m},D^*_{4n}/C_{2n})$ & $4mn$ &   \\ 
 4. & $(C_{4m}/C_{2m},D^*_{8n}/D^*_{4n})$ & $8mn$ &   \\ 
 5. & $(C_{2m}/C_{2m},T^*/T^*)$ & $24m$ &   \\
 6. & $(C_{6m}/C_{2m},T^*/D^*_{8})$ & $24m$ &   \\ 
 7. & $(C_{2m}/C_{2m},O^*/O^*)$ & $48m$ &   \\
 8. & $(C_{4m}/C_{2m},O^*/T^*)$ & $48m$ &   \\ 
 9. & $(C_{2m}/C_{2m},I^*/I^*)$ & $120m$ &   \\ 
 10. & $(D^*_{4m}/D^*_{4m},D^*_{4n}/D^*_{4n})$ & $8mn$ &   \\
 11. & $(D^*_{4mr}/C_{2m},D^*_{4nr}/C_{2n})_s$ & $4mnr$ & $\gcd(s,r)=1$ \\ 
 $11^{\prime}$. & $(D^*_{2mr}/C_{m},D^*_{2nr}/C_{n})_s$ & $mnr$ &  $\gcd(s,r)=1$ $\gcd(2,n)=1$  \\
&&&  $\gcd(2,m)=1$ $\gcd(2,r)=2$\\
 12. &  $(D^*_{8m}/D^*_{4m},D^*_{8n}/D^*_{4n})$ & $16mn$ & \\
 13. &  $(D^*_{8m}/D^*_{4m},D^*_{4n}/C_{2n})$ & $8mn$ & \\
 14. &  $(D^*_{4m}/D^*_{4m},T^*/T^*)$ & $48m$ & \\
 15. &  $(D^*_{4m}/D^*_{4m},O^*/O^*)$ & $96m$ & \\
16. &  $(D^*_{4m}/C_{2m},O^*/T^*)$ & $48m$ & \\
17. &  $(D^*_{8m}/D^*_{4m},O^*/T^*)$ & $96m$ & \\
18. & $(D^*_{12m}/C_{2m},O^*/D^*_{8})$ & $48m$ & \\
19. & $(D^*_{4m}/D^*_{4m},I^*/I^*)$ & $240m$ & \\
20. & $(T^*/T^*,T^*/T^*)$ & $288$ & \\
21. & $(T^*/C_2,T^*/C_2)$ & $24$ & \\
$21^{\prime}.$ & $(T^*/C_1,T^*/C_1)$ & $12$ & \\
22. & $(T^*/D^*_{8},T^*/D^*_{8})$ & $96$ & \\
23. & $(T^*/T^*,O^*/O^*)$ & $576$ & \\
24. & $(T^*/T^*,I^*/I^*)$ & $1440$ & \\
25. & $(O^*/O^*,O^*/O^*)$ & $1152$ & \\
26. & $(O^*/C_2,O^*/C_2)$ & $48$ & \\
$26^{\prime}.$ & $(O^*/C_1,O^*/C_1)_{Id}$ & $24$ & \\
$26^{\prime\prime}.$ & $(O^*/C_1,O^*/C_1)_f$ & $24$ & \\
27. & $(O^*/D^*_{8},O^*/D^*_{8})$ & $192$ & \\
28. & $(O^*/T^*,O^*/T^*)$ & $576$ & \\
29. & $(O^*/O^*,I^*/I^*)$ & $2880$ & \\
30. &   $(I^*/I^*,I^*/I^*)$ & $7200$ & \\
31. &   $(I^*/C_2,I^*/C_2)_{Id}$ & $120$ & \\
$31^{\prime}.$ & $(I^*/C_1,I^*/C_1)_{Id}$ & $60$ & \\
32. &   $(I^*/C_2,I^*/C_2)_{f}$ & $120$ & \\
$32^{\prime}.$ & $(I^*/C_1,I^*/C_1)_f$ & $60$ & \\
33. & $(D^*_{8m}/C_{2m},D^*_{8n}/C_{2n})_f$ & $8mn$ &   $m\neq 1$  $n\neq 1$. \\ 
 $33^{\prime}$. & $(D^*_{8m}/C_{m},D^*_{8n}/C_{n})_f$ & $4mn$ &  $\gcd(2,n)=1 \gcd(2,m)=1$ \\
&&& $m\neq 1$ and $n\neq 1$.   \\
34. &  $(C_{4m}/C_{m},D^*_{4n}/C_{n})$ & $2mn$ & $\gcd(2,n)=1 \gcd(2,m)=1$ \\
\hline
\end{tabular}

\end{table}

Note that, in most cases, the group is entirely determined by $(L, L_K, R, R_K)$ up to conjugacy, as the isomorphism $\phi$ is unique or does not change the conjugacy class, so we use Du Val's notation $(L/ L_K, R/ R_K)$, sometimes using a subscript to denote the isomorphism $\phi$.

In case of Family 1, we write $(C_{2mr}/ C_{2r}, C_{2nr}/C_{2r})_s$ (with $\gcd(s,r) = 1$) to indicate that the isomorphism $\phi : L/L_K \rightarrow R/R_K$ is the one defined by $\phi(e^{2i\pi/2mr}C_{2m})= e^{2i\pi s/2nr}C_{2n}$. Similarly, in case of Family $1'$, the group  $(C_{mr}/C_m, C_{nr}/C_n)_s$ (with $m,n$ odd, $r$ even and $\gcd(s,r)=1$) is given by  the isomorphism $\phi : L/L_K \rightarrow R/R_K$ such that  $\phi(e^{2i\pi/2mr}C_{m})= e^{2i\pi s/nr}C_{n}$.

In case of Family 11, the isomorphism related to $(D^*_{4mr}/C_{2m}, D^*_{4nr}/C_{2n})_s$ is determined by $\phi(e^{2i\pi/2mr}C_{2m})= e^{2i\pi s/2nr}C_{2n}$ and $\phi(j)=j$.  In case of Family $11'$, the isomorphism related to $(D^*_{2mr}/C_{m}, D^*_{2nr}/C_{n})_s$ is determined by $\phi(e^{2i\pi/mr}C_{m})= e^{2i\pi s/2nr}C_{n}$ and $\phi(j)=j$. For Families $26^\prime,\, 26^{\prime\prime},\,31,\, 31^\prime,\,32,\, 32^\prime,\,33$ and  $33^\prime$, the definition of the isomorphism $\phi$ can be found in \cite{mecchia-seppi2}. If we add the word ``bis'' to the number of the family we are referring to the family when the roles of $L$ and $R$ are switched with respect to the original one. A group is conjugate to the corresponding one in the ``bis'' family by an orientation reversing isometry, but in $\SO(4)$ the two groups might be non-conjugate. More details can be found in \cite[Section 3.2]{mecchia-seppi2}.

\medskip

Most of the quotients of $S^3$ by one of these groups are not manifolds so we need to discriminate between the groups acting freely on $S^3$, whose associated quotient is a manifold, and the others, whose associated quotient is an orbifold with nonempty singular set. The subgroups of $\SO(4)$ acting freely have been classified by Seifert and Threlfall in \cite{threlfall-seifert,threlfall-seifert2}. A description of these groups is present in \cite{mccullough} where the notation is similar to ours except for Families 1 and $1'.$ For this case we use the results of \cite{mecchia-seppi} to determine when the singular locus is empty.

Finally, we find that the subgroups of $SO(4)$ acting freely are those of 

\begin{itemize}
\item Family 1, $\tilde{G}= (C_{2mr}/C_{2m},C_{2nr}/C_{2n})_s$ with $\gcd(r,s) =1$, $\gcd(m,n)=1$, $mn$ even or $r$ odd, and $\gcd(n-sm,2mnr)= \gcd(n+sm,2mnr)$
\item Family $1'$, $\tilde{G}= (C_{mr}/C_{m},C_{nr}/C_{n})_s$ with $m,n$ odd, $r$ even, $\gcd(r,s)=1$, $\gcd(m, n)=1$, and $\gcd(n-sm,mnr)= \gcd(n+sm,mnr)$
\item Family 2, $\tilde{G}= (C_{2m}/C_{2m},D^*_{4n}/D^*_{4n})$ with $\gcd(m,2n) = 1$
\item Family 3, $\tilde{G}= (C_{4m}/C_{2m},D^*_{4n}/C_{2n})$, with $\gcd(m,n)= 1$ and $m$ even
\item Family 5, $\tilde{G}= (C_{2m}/C_{2m},T^*/T^*)$, with $\gcd(m, 6) = 1$
\item Family 6, $\tilde{G}= (C_{6m}/C_{2m}, O^*/T^*)$, with $m$ odd, $3|m$
\item Family 7, $\tilde{G}= (C_{2m}/C_{2m},O^*/O^*)$, with $\gcd(m, 6) = 1$
\item Family 9, $\tilde{G}= (C_{2m}/C_{2m},I^*/I^*)$, with $\gcd(m, 30) = 1$
\end{itemize}

In accordance to the family of $\tilde{G}$, the manifold $S^3/G$ is called a \textit{prism manifold}  if $\tilde{G}$ is in Family 2 or 3, a \textit{tetrahedral manifold} if $\tilde{G}$ is in Family 5 or 6, an \textit{octahedral manifold} if $\tilde{G}$ is in Family 7, and an \textit{icosahedral manifold} if $\tilde{G}$ is in Family 9. If $\tilde{G}$ is in Family 1 or $1'$, then $S^3/G$ is a lens space.


In the following, the study of those groups will in most cases be separated between the general case and the cases of "small indices", that is when one of the groups involved is $C_2$ or $D^*_8$, or when $r\in \{1,2\}$ for Families 1 and $1'$.

\section{Isometry groups} \label{isometry groups}

An isometry of the spherical 3-orbifold $O= S^3/G$  can be lifted to an isometry of $S^3$  normalizing $G$; if the initial isometry is orientation-preserving, then the lift to $S^3$ is orientation-preserving too. This implies that $\Isompr(S^3/G)\cong \Norm_{\SO(4)}(G)/G$. By the 2:1 correspondence $\Phi$ used in the previous section, we can deduce that $\Isompr(S^3/G)\cong \Norm_{S^3\times S^3}(\tilde{G})/\tilde{G}$, 

Thus, the isometry groups of spherical orbifolds can be computed by analyzing the normalizer of the finite subgroups of $S^3\times S^3$. In \cite{mccullough} the case of 3-manifolds is considered, while the isometry groups of all spherical orbifolds are computed in \cite{mecchia-seppi2}. 

Let us recall that the finite subgroups of $S^3$ are known and that their normalizers are:
\begin{align*}
 & \Norm_{S^3}(C_n)=O(2)^* \text{ if } n>2, \text{ and }  \Norm_{S^3}(C_2)=S^3\\
 & \Norm_{S^3}(D^*_{4n})=D^*_{8n} \text{ if } n>2, \text{ and } \Norm_{S^3}(D^*_8)=O^*\\
 & \Norm_{S^3}(T^*)=O^* \\
 & \Norm_{S^3}(O^*)=O^* \\
 & \Norm_{S^3}(I^*)=I^*
\end{align*}

The group $O(2)^*$ is generated by the subgroup $S^1=\{e^{i\theta},\,\theta \in \mathbb{R}\}$ and $j;$ this group is a central extension of $\O(2)$ by a subgroup of order two.

The computation of the normalizer of a group $\tilde{G}= (L,L_K, R, R_K, \phi)$ is given by the following lemma that is a special case of Lemma~\ref{conjugation}.

\begin{lemma}
Let $\tilde{G}=(L, L_K , R, R_K , \phi)$ be a finite subgroup of $S^3\times S^3$ containing $(-1,-1)$. An element $(g, f )\in S^3\times S^3$ normalizes $\tilde{G}$ if and only if the following conditions are satisfied:
\begin{enumerate}
\item $g \in \Norm_{S^3}(L)\cap \Norm_{S^3}(L_K)$
\item $f \in \Norm_{S^3}(R)\cap \Norm_{S^3}(R_K)$
\item The equality $\bar{c}_f \circ \phi = \phi \circ \bar{c}_g$ holds, where $\bar{c}_g$ (resp. $\bar{c}_f$) is the application $L/L_K \rightarrow L/L_K$ (resp. $R/R_K \rightarrow R/R_K$) induced in the quotient by  $c_g: L\rightarrow L$ (resp. $c_f: R\rightarrow R$), the conjugation by $g$ (resp. by $f$).
\end{enumerate}
\end{lemma}

This lemma can be used to compute the orientation-preserving isometry groups of the spherical 3-manifolds and 3-orbifolds. The isometry groups of the manifolds are reported in  Table~\ref{isometry groups} where we use the same notation of \cite{mecchia-seppi2}. We must distinguish between the general case and some cases involving small indices. In the general case, we suppose that $m,n>1$ and $r>2.$ The symbol $\widetilde{\times}$ indicates a central product where the two central involutions of the factors are identified. 

\begin{table}
\caption{Isometry groups}
\label{isommetry-groups}
\begin{tabular}{|l|c|c||l|c|c|}
\hline
 &  $\tilde G$ (general case) & $\Isompr(S^3/G)$ &  & $\tilde{G}$ (small indices) & $\Isompr(S^3/G)$\\
\hline
 1. & $(C_{2mr}/C_{2m},C_{2nr}/C_{2n})_s$ & ${Dih}(S^1\times S^1)$  
 & 1. & $(C_{2m}/C_{2m},C_{2n}/C_{2n})$ & $\O(2)\times \O(2)$  \\  
 $1^{\prime}$. & $(C_{mr}/C_{m},C_{nr}/C_{n})_s$ & ${Dih}(S^1\times S^1)$
 & &  $(C_{4m}/C_{2m},C_{4n}/C_{2n})$ & $\O(2)\widetilde{\times}\O(2)$ \\
 2. & $(C_{2m}/C_{2m},D^*_{4n}/D^*_{4n})$ & $\O(2)\times \Z_2$   
 & & $(C_{2}/C_{2},C_{2n}/C_{2n})$ & $\SO(3)\times \O(2)$ \\
 3. & $(C_{4m}/C_{2m},D^*_{4n}/C_{2n})$ & $\O(2)\times \Z_2$ 
 & & $(C_2/C_2,C_2/C_2)$ & $\mathrm{P}\SO(4)$\\
 5. & $(C_{2m}/C_{2m},T^*/T^*)$ & $\O(2)\times\Z_2$
 & $1^{\prime}$. & $(C_{2m}/C_{m},C_{2n}/C_{n})$ & $\O(2)^*\widetilde{\times} \O(2)^*$ \\ 
 6. & $(C_{6m}/C_{2m},T^*/D^*_{8})$ & $\O(2)$  
 & &  $(C_2/C_1,C_{2n}/C_n)$ & $S^3\widetilde{\times}\O(2)^*$ \\
 7. & $(C_{2m}/C_{2m},O^*/O^*)$ & $\O(2)$  
 & 2. & $(C_{2}/C_{2},D^*_{4n}/D^*_{4n})$ & $\SO(3)\times \Z_2$  \\
 9. & $(C_{2m}/C_{2m},I^*/I^*)$  &  $\O(2)$
 &  &  $(C_{2m}/C_{2m},D^*_{8}/D^*_{8})$   & $\O(2)\times D_6$  \\ 
 &  &  
 &  &  $(C_{2}/C_{2},D^*_{8}/D^*_{8})$ &  $\SO(3)\times D_6$ \\
 &  &  
 & 5. & $(C_{2}/C_{2},T^*/T^*)$  & $\SO(3)\times\Z_2$ \\
 &  &  
 &  7. &  $(C_{2}/C_{2},O^*/O^*)$ &  $\SO(3)$\\
 &  &  
 &9. & $(C_{2}/C_{2},I^*/I^*)$ & $\SO(3)$\\
 \cline{1-6}
 \end{tabular}
\end{table}

\section{Involutions of spherical manifolds}\label{sec: involutions}

As indicated in Table~\ref{isommetry-groups}, the orientation-preserving isometry groups of spherical manifolds can be expressed using the groups $\mathbb{Z}_2,$ $\mathbb{D}_6,$ $\O(2),$ $\O(2)^*,$ $ S^3,$ $ \SO(3),$ $ \SO(4),$ $ \PSO(4)$ and $\Dih(S^1\times S^1)$,  being either one of those or a product (direct or central) of two of them.  It is useful to notice that $O(2)$ is a dihedral extension of $S^1$.

If $G$ is a group, we denote by $\In(G)$ the set of involutions of $G$, and by $\I(G)$ the set of conjugacy classes of the involutions. In Table~\ref{involutions}, we describe a complete set of representatives of $\I(G)$ for each of these groups $G$, that is a set containing exactly one element for each conjugacy class of involutions. We recall the two conjugate involutions give diffeomorphic quotients, thus in this context we are interested in the classification of involutions up to conjugacy.
\\
In Table~\ref{involutions}, we denote by $c$ an involution acting dihedrally on the normal subgroup $S^1$ in the group $\O(2)$, and we use the following notations for the groups  $ \SO(3),$ $ \SO(4)$ and $ \PSO(4)$:

\[
  R_\theta =
  \left[ {\begin{array}{ccc}
    1 & 0 & 0 \\
    0 & \cos(\theta) & -\sin(\theta) \\
    0 & \sin(\theta) & \cos (\theta) \\
  \end{array} } \right]
  = \left[ {\begin{array}{cc}
    1 & 0_{1\times2}  \\
    0_{2\times 1} & r_\theta\\
  \end{array} } \right] 
  \text{ with } r_\theta =
  \left[ {\begin{array}{cc}
    \cos \theta & -\sin \theta  \\
    \sin \theta & \cos \theta\\
  \end{array} } \right] 
\]

\[
  T_{\theta,\rho} =
  \left[ {\begin{array}{cc}
    r_\theta & 0  \\
    0 & r_\rho \\
  \end{array} } \right]
\]


We now analyse the cases when the isometry group is the product of two groups. It is easy to see that if $F,G$ are two groups, then $\In(F\times G) = \In(F) \times \In(G)\cup \{1_F\} \times \In(G) \cup \In(F) \times \{1_G\}$. 
Furthermore, in a direct product two involutions $(a,b)$ and $(a',b')$ are conjugate if and only if $a$ and $b$ are conjugate to  $a'$ and $b'$, respectively.  If $\mathcal{I}(F)$ (resp. $\mathcal{I}(G)$) is a complete set of representatives of $\I(F)$ (resp. $\I(G)$), then the union of   $ \mathcal{I}(F) \times \mathcal{I}(G),$ $ \{1_F\} \times \mathcal{I}(G)$ and $\mathcal{I}(F) \times \{1_G\}$ is a complete set of representatives of $\I(F\times G)$. These remarks can be used to compute the involutions and their conjugacy classes in the isometry groups that are direct products of groups in Table~\ref{involutions}.

Consider now $A$ and $B$, two groups containing an involution in the centre. We use $-1$ to denote the central involution of both $A$ and $B$. We want to compute the involutions and their conjugacy classes of the central product $A\widetilde{\times}B = A\times B / \{ \pm (1,1)\}$. The central products appearing in Table~\ref{isommetry-groups} are of this type.

We can see that an element $(x, y)$ in $A\times B$, gives an involution in the quotient $A\widetilde{\times}B$ if and only if $(x,y)$ belongs to $\In(A\times B) \setminus \{(-1,-1)\}$  or  to $\Rac(A)\times \Rac(B)$, where $\Rac(G) = \{z\in G| z^2=-1\}$. We note also that elements  in  $\In(A\times B) \setminus \{(-1,-1)\}$ can not be conjugate to elements of $Rac(A)\times Rac(B)$. We denote by $R(G)$ the conjugacy classes of the elements in $\Rac(G).$

In the last column of Table~\ref{involutions}, we provide a complete set of representatives for $R(G)$ for the groups involved as factors in central products.

The central products appearing as isometry groups are: $\O(2)\widetilde{\times}\O(2),$ $\O(2)\widetilde{\times}\O(2)$ and $S^3\widetilde{\times}\O(2)^*.$ We consider as example the group   $\O(2)\widetilde{\times}\O(2).$ In the group $\O(2)\times\O(2)$ we have 8 involutions up to conjugacy: $(-1,-1)$, $(-1,c),$ $(c,-1),$ $(c,c),$ $(1,-1),$ $(1,c),$ $(c,1),$ and  $(-1,1).$
In the quotient $\O(2)\widetilde{\times}\O(2)= \O(2)\times\O(2) / \{ \pm (1,1)\}$ we have that:
\begin{itemize}
\item the element  $(-1,-1)$ gives the trivial element;
\item the cosets of $(-1,1)$ and $(1,-1)$ coincide;
\item the cosets of $(-1,c)$ and $(1,-c)$ coincide and this implies that the cosets of  $(-1,c)$ and $(1,c)$ are conjugate;  
\item analogously the cosets of $(c,-1)$ and $(c,1)$ are conjugate.
\end{itemize}

Finally, in $\O(2)\widetilde{\times}\O(2)$ a complete set of representatives of the set of the conjugacy classes of involutions is  given by the cosets of the following elements:
$(-1,1)$, $(1,c)$, $(c,1)$, $(c,c)$, $(i,i).$

An analogous computation can be made for the other central products. For $\O(2)^*\widetilde{\times}\O(2)^*$ a complete set of representatives is given by the cosets of $(1,-1),$ $(i,i),$ $(j,j),$ $(i,j),$ and $(j,i).$

For the group $S^3\widetilde{\times}\O(2)^*$, we have three conjugacy classes of involutions given by the cosets of $(1,-1),$ $(i,i)$ and  $(i,j).$

The conjugacy classes of involutions are summarised in Tables~\ref{general quotient} and \ref{small index quotient}: each line of the second column corresponds to a conjugacy class of involutions. The only exception is Family 1 in Table~\ref{general quotient}, which has four conjugacy classes of involutions, but three are treated together. These tables show the geometric features of the involutions. How to read and complete them is explained in the following sections.

\begin{table} 
\caption{Conjugacy classes of involutions}
\label{involutions}
\begin{tabular}{|c|c|c|c|}
\hline
group $G$ & representation  &  \multicolumn{2}{|c|}{complete set of representatives of: } \\[1ex]
 & of the group  $G$  &  $I(G)$ & $R(G)$ \\[1ex]
\hline
\hline
&&&\\
$\mathbb{Z}_2$ & $\{[0],[1]\}$ & [1] &\\[1ex]
$\mathbb{D}_6$ & $\{1,r,r^2, s, sr, sr^2\}$ & $s$ &\\[1ex]
$S^1$ & $\{e^{i\theta}, \theta \in \mathbb{R}\}$ & $-1$ & $i$ \\[1ex]
$O(2)$ & $S^1\cup cS^1$  & $-1,c$ & $i$\\[1ex]
$O^*(2)$ & $S^1\cup S^1j$ & $-1$ & $i,j$\\[1ex]
$S^3$ & $\{q\in \mathbb{H} | \|q\|=1\}$ & $-1$ & $i$ \\[1ex]
$SO(3)$ & & $R_\pi$ & \\[1ex]
$SO(4)$ & & $T_{0,\pi}, T_{\pi, 0}, T_{\pi,\pi}$ &\\[1ex]
$PSO(4)$ & $\{[Q], Q\in SO(4)\}$ & $ [T_{0,\pi}], [T_{\pi, 0}], [T_{\pi/2, -\pi/2}], [T_{\pi/2, \pi/2}]$  &\\[1ex]
$Dih(S^1\times S^1)$  & $(S^1\times S^1) \rtimes \mathbb{Z}_2$ & $((-1,-1),[0]),((-1,1),[0]),$ & \\[1ex]
&  & $((1,-1),[0]), ((1,1),[1])$ & \\[1ex]

\hline
\end{tabular}
\end{table}

\section{Geometric properties of involutions}

The advantage of this approach is the possibility of easily obtaining geometric information about involutions; for example,  we can identify which act freely and which are hyperelliptic. Most of the geometric properties can be deduced by the quotient of the manifold by the involution.

\subsection{Extension given by involutions.} \label{extensions}

Given a closed spherical 3-manifold $S^3/G$ and an involution $f\in \Isompr(S^3/G)$, we have $(S^3/G)/\langle f \rangle = S^3/\langle G, \bar{f} \rangle$, where $\bar{f}$ is any lift of $f$ to $\Norm_{\SO(4)}(G).$ 
Any quotient orbifold $S^3/\langle G, \bar{f}\rangle$ admits a Seifert fibration and the formulae in \cite{mecchia-seppi} give the Seifert invariants. To use these formulae we work  in $S^3\times S^3$; we note that $\langle G, \bar{f} \rangle= \Phi (\langle\tilde{G}, \tilde{f} \rangle)$, where $\tilde{f}$ is any lift of $f$ to $\Norm_{S^3\times S^3}(\tilde{G})$, i.e. a preimage of $\bar{f}$ under  $\Phi.$  

We recall that the conjugacy classes of involutions can be computed following the strategy presented in Section~\ref{sec: involutions}.

In this section, we will try to be consistent with the following notations: $G$ is a finite subgroup of $\SO(4)$, $\tilde{G}$ is the associated subgroup of $S^3\times S^3$, $f$ is an involution of $\Isompr(S^3/G)$, $\bar{f} $ is one of its lifts to $\Norm_{\SO(4)}(G),$ and  $\tilde{f}$ is one of the lifts of $f$ to $\Norm_{S^3\times S^3}(\tilde{G})$; for the sake of brevity, we often denote $\langle\tilde{G}, \tilde{f} \rangle$ by $H$.


\begin{remark}
\label{computations-about-extensions}

We know that  $\langle\tilde{G}, \tilde{f}\rangle = \tilde{G} \cup \tilde{f}\tilde{G}$ is an extension of index 2 of $\tilde{G}$, but also that if $\tilde{G}=(L, L_K, R, R_K, \phi)$, $\langle\tilde{G}, \tilde{f}\rangle = (L', L_K', R', R_K', \phi')$ and $\tilde{f}=(p,q)$, then $L'=L\cup pL$ and $R'=R\cup qR$. 
\\
Moreover, we can also compute $L_K'$ and $R_K'$: we denote by $\pi_R:R\rightarrow R/R_K$ (resp. $\pi_L:L\rightarrow L/L_K$) the projection from $R$ (resp. $L$) to the quotient $R/R_K$ (resp. $L/L_K$).

\begin{itemize}
\item If $p\in L$, then an element $z$ of $S^3$ belongs to $R_K'$ if and only if $z$ belongs to $R_K$ or there exists $ (x,y)\in \tilde{G}$ such that $ (p,q)(x,y)=(1,z)$. We obtain the following equivalences:

\begin{center}
\begin{tabular}{lcl}
    $\exists (x,y)\,\in \tilde{G}$ s.t. $(p,q)(x,y)=(1,z) $    &
     $\Leftrightarrow$  &
      $\exists (x,y)\in \tilde{G}$  s.t.  $x=p^{-1}$ and $qy=z$ \\
       &  $\Leftrightarrow$  &
     $\exists y\in R$ s.t. $(p^{-1},y)\in \tilde{G}$ and $qy=z$ \\
      &  $\Leftrightarrow$  &
      $\exists  y \in R $ s.t. $\pi_R(y)= \phi(\pi_L(p^{-1}))$ and $qy=z$ \\
      &  $\Leftrightarrow$  &
       $\exists  y \in \pi_R^{-1}(\phi(\pi_L(p^{-1}))) $ s.t. $z=qy.$ 
\end{tabular}       
\end{center}
\smallskip

    Therefore we have that $R_K'=R_K \cup q\,\pi_R^{-1}\circ\phi\circ\pi_L(p^{-1})$;  moreover the union is disjoint and $|R_K'|=2|R_K|$ holds, otherwise there would be some $y\in \pi_R^{-1}(\phi(\pi_L(p^{-1})))$ such that $qy$ belongs to $R_K$, and then we would have $(p^{-1},y) \in \tilde{G}$ and $(1, qy) \in \tilde{G}$, so $(p,q)= (1, qy) (p^{-1},y)^{-1} \in \tilde{G}$, which is impossible.
Similarly, if $q\in R$, then $L_K'= L_K \cup p\,\pi_L^{-1}(\phi^{-1}(\pi_R(q^{-1})))$ and $|L_K'|=2|L_K|$.

\item Conversely, if $p\notin L$ (resp. $q\notin R$), then $|L'|=2|L|$ (resp. $|R'|=2|R|$), and there is no element $(x,y)$ in $\tilde{G}$ such that $x=p^{-1}$, so $R_K'=R_K$ (resp. $L_K'=L_K$).
\end{itemize}

\end{remark}

Thus, we obtained a complete description (safe for the isomorphism $\phi'$) of $H=\langle\tilde{G}, \tilde{f}\rangle= (L', L_K', R',$ $R_K', \phi')$. The other data determine the isomorphism except for Families $1$ and $1'$. We start describing the computation for  Family 2 when $m=1$.  Then we consider the most complicated cases of Families 1 and $1'$ in detail. The other families can be treated by following the same strategy. 

The results of these computations are contained in Tables~\ref{fibrations-one}, \ref{fibrations-two} and \ref{fibrations-three}. For each group $G$, the tables give an extension  $\langle\tilde{G}, \tilde{f}\rangle$ for each conjugacy class of involutions.   For $a,b, r \in \mathbb{Z}$ we write $a|b$ if $a$ divides $b$  and $a \equiv b [r]$ if $a$ and $b$ are congruent modulo $r$. The remainder of the euclidean division of $a$ by $b$ is denoted by $a\%b$. The tables considering small indices contain some extensions that are in Families $1'$ and $1'$ with $r=4$; in this case, the subscript $s$ can be omitted since we can replace $s$ with $r-s$ obtaining  a conjugate group (see~\cite[page 803]{mecchia-seppi})
 
\bigskip
\noindent
\textbf{Family 2:} $\tilde{G}=(C_2/C_2,D^*_{4n}/D^*_{4n}).$
We have $\Norm_{S^3\times S^3}(\tilde{G})= (S^3/S^3,D^*_{8n}/D^*_{8n})$, and $\Isompr(S^3/G)=\Norm_{S^3\times S^3}(\tilde{G}) /\tilde{G}= SO(3)\times \mathbb{Z}_2$. Therefore, $\Isompr(S^3/G)$ has three conjugacy classes of involutions, whose representatives are $(R_\pi, 0), (I_3, 1)$ and $(R_\pi, 1)$.
\begin{itemize}
\item If $f= (R_\pi, 0)$ then $\tilde{f}=(i,1)$ is a lift of $f$ to $S^3\times S^3$. Therefore we have $H= \langle (C_2/C_2,D^*_{4n}/D^*_{4n}), (i,1)\rangle= (C_4/C_4, D^*_{4n}/D^*_{4n}).$ 
\item If $f= (I_3,1)$ then $\tilde{f} = ( 1 , e^{i\pi/4n} )$ is a lift of $f$ to $S^3\times S^3$.  Therefore we have $H= ((C_2/C_2,D^*_{4n}/D^*_{4n}) , ( 1 , e^{2i\pi/4n} )\rangle = (C_2/C_2, D^*_{8n}/D^*_{8n})$
\item Of $f= (R_\pi, 1)$ then $\tilde{f}=(i,e^{i\pi/4n})$ is a lift of $f$ to $S^3\times S^3$. Therefore we have $H = \langle (C_2/C_2,D^*_{4n}/D^*_{4n}) ,( i , e^{2i\pi/4n} ) ) \rangle =  ( C_4 / C_2, D^*_{8n} / D^*_{4n} ),$ because $i\notin C_2$ and $e^{2i\pi/4n} \notin D^*_{4n}$.

\end{itemize}

\bigskip

\noindent
\textbf{Family 1$'$:} $\tilde{G}= (C_{mr}/C_{m}, C_{nr}/C_{n})_s.$ We recall that $S^3/G$ is a manifold if and only if    $\gcd(r,s)=1,$ $r$ even, $m$ odd, $n$ odd, $\gcd(m,n)=1$ and $\gcd(n-sm,mnr)= \gcd(n+sm,mnr).$ 
We have $\Norm_{S^3\times S^3}(\tilde{G})= (O(2)^*/S^1, O(2)^*/S^1)$ and $\Isompr(S^3/G)= \Norm_{S^3\times S^3}(\tilde{G})/G= Dih(S^1\times S^1/ \tilde{G}) \cong Dih(S^1\times S^1)$, seen as $(S^1\times S^1) \rtimes \mathbb{Z}_2$. The conjugacy classes of the involutions are represented by $((-1,1),0), ((1,-1,), 0), ((-1,-1),0)$ and $((1,1),1)$. As explained in \cite{mecchia-seppi2}, the isomorphism between $S^1\times S^1$ and $S^1\times S^1/\tilde{G}$ is given by

$$\gamma :S^1\times S^1 \rightarrow S^1\times S^1/\tilde{G} \quad (e^{i\alpha}, e^{i\beta}) \rightarrow \left(e^{i\left(\frac{\alpha}{mr}+ \frac{\beta}{m}\right)}, e^{i\left(\frac{s\alpha}{nr}+\frac{(s+1)\beta}{n}\right)}\right)\tilde{G}$$
Therefore, we have:
$$\gamma(-1,1) = \left(e^{i\frac{\pi}{mr}}, e^{i\frac{s\pi}{nr}}\right)\tilde{G},\, 
\gamma(1,-1)=\left(e^{i\frac{\pi}{m}}, e^{i\frac{(s+1)\pi}{n}}\right)\tilde{G}\text{ and}$$  
$$ \gamma(-1,-1)=\left(e^{i\frac{\pi}{mr}(1+r)}, e ^{i\frac{\pi}{nr}(rs+r+s)}\right)\tilde{G}$$ 
which give us the lifts of $((-1,1),0), ((1,-1),0)$ and $((-1,-1),0)$ respectively. The lift of $((1,1),1)$ is $(j,j)$.
\begin{itemize}
\item Suppose $\tilde{f}=(e^{i\pi/mr}, e^{is\pi/nr})$. We denote by $(L'/L_K', R'/R_K')_{s'}$ the representation of  $H=\tilde{G}\cup \tilde{f}\tilde{G}$ as subgroup of $S^3\times S^3$. Since $e^{i\pi/mr} \in C_{2mr} \setminus C_{mr}$, we have $L'= C_{2mr}$. For the same reason, we also have $R'= C_{2nr}$. By the discussion in Remark~\ref{computations-about-extensions}, we have $L_K'=L_K=C_m$ and $R_K'=R_K=C_n$
In order to find $s'$, it is sufficient to find an element $y$ such that $(e^{i\pi/mr},y) \in H$, and we note that $y= e^{is\pi/nr}$ works, so $s'=s$. We obtain groups in Family $1'$.

\item Suppose $\tilde{f}= (e^{i\pi/m}, e^{i(s+1)\pi/n})$. With the same notations as above, we have $e^{i\pi/m}= e^{2i\pi (r/2)/mr} \in C_{mr}$ because $r$ is even, and therefore $L'=C_{mr}$. Symmetrically, $R'=C_{nr}$.
To find $L'_K$ and $R'_K$ we can either use Remark~\ref{computations-about-extensions} or compute it directly: if $z\in R_K'$, then either $z\in R_K$ or there exists $(x,y)\in \tilde{G}$ such that  $\tilde{f}\cdot(x,y)=(1,z)$, and then $x=e^{-i\pi/m}$, so $y\in e^{-is\pi/n}C_n$, hence $z\in e^{i\pi/n}C_n$ ($\subset C_{2n}$). Conversely, $\{1\}\times e^{i\pi/n}C_n$ is contained in $\tilde{f}\tilde{G}$, which is why $R_K'= C_{2n}$. Similarly, $L_K'=C_{2m}$.
Finally, we have $s'=s \% (r/2)$ since $(e^{2i\pi/2m(r/2)}, e^{2i\pi s/2n(r/2)})\in H$, so $H = (C_{2m(r/2)} / C_{2m},
  C_{2n(r/2)} /  C_{2n})_{s\%(r/2)}.$ These groups belong to Family 1.

\item Suppose $\tilde{f}=(e^{i\pi(1+r)/mr}, e^{i\pi (rs+r+s)/nr})$. With the same reasoning as in the first case, we find $L'=C_{2mr}, R'=C_{2nr}, L_K'=C_m$ and $R_K'=C_n$.
To find $s'$, we search $y$ such that $(e^{i\pi/mr}, y) \in H$. We have $(e^{2i\pi /mr}, e^{2i\pi s/nr})\in H$, so 
$$(e^{i\pi(1+r)/mr}, e^{i\pi (rs+r+s)/nr}) (e^{2i\pi /mr}, e^{2i\pi s/nr})^{-r'}=(e^{2i\pi /2mr}, e^{2i\pi (r+s)/2nr})\in H,$$ and then $s'=r+s$.
Finally, we get  $H = (C_{2mr}/C_m, C_{2nr}/ C_n)_{s+r}$ which are groups in Family $1'$.

\item If $\tilde{f}= (j,j)$, then $H= (D_{2mr}^*/C_{m}, D^*_{2nr}/C_{n})_s$ and we are in Family $11'.$
\end{itemize}
\bigskip

\noindent
\textbf{Family 1:} Consider $\tilde{G}= (C_{2mr}/C_{2m}, C_{2nr}/C_{2n})_s$, with  $\gcd(r,s)=1,$ $\gcd(m,n)=1$, $mn$ even or $r$ odd, and $\gcd(n-sm,2mnr)= \gcd(n+sm,2mnr).$ 

We have $\Norm_{S^3\times S^3}(\tilde{G})= (O(2)^*/S^1, O(2)^*/S^1)$ and $\Isompr(S^3/G)= \Norm_{S^3\times S^3}(\tilde{G})/G= Dih(S^1\times S^1/ \tilde{G})= Dih(S^1\times S^1)$. The situation is similar to Family $1'$, and the conjugacy classes are the same.
The isomorphism between $S^1\times S^1$ and $S^1\times S^1/\tilde{G}$ is given by

$$\gamma :S^1\times S^1 \rightarrow S^1\times S^1/\tilde{G}\quad  (e^{i\alpha}, e^{i\beta}) \rightarrow \left(e^{i\left(\frac{\alpha}{2mr}+ \frac{\beta}{2m}\right)}, e^{i\left(\frac{s\alpha}{2nr}+\frac{(s+1)\beta}{2n}\right)}\right)\tilde{G}$$
Therefore, we have: 
$$\gamma(-1,1) = (e^{i\pi/2mr}, e^{is\pi/2nr})\tilde{G},\,\gamma(1,-1) = (e^{i\pi/2m}, e^{i(s+1)\pi/2n})\tilde{G} \text{ and }$$ $$\gamma(-1,-1) = (e^{i\pi (1+r)/2mr}, e^{i\pi(rs+r+s)/2nr})\tilde{G},$$ which give us the lifts of $((-1,1),0), ((1,-1),0)$ and $((-1,-1),0)$, respectively. The lift of $((1,1),1)$ is $(j,j)$.

Let us first consider what happens when $r$ is even. Everything works the same way as in the case of Family $1'$, and so we get:
\begin{itemize}
\item If $\tilde{f}= (e^{i\pi/2mr}, e^{is\pi/2nr})$ then $H= (C_{2m\cdot 2r}/C_{2m}, C_{2n\cdot 2r}/C_{2n})_s$ (Family 1).  
\item If $\tilde{f}= (e^{i\pi/2m}, e^{i(s+1)\pi/2n})$ then  $H= (C_{4m(r/2)}/C_{4m}, C_{4n(r/2)}/C_{4n})_{s\%(r/2)}$ (Family 1). 
\item If $\tilde{f}= (e^{i\pi (1+r)/2mr}, e^{i\pi (rs+r+s)/2nr})$ then $H= (C_{4mr}/C_{2m}, C_{4nr}/C_{2n})_{r+s}$ (Family 1). 
\item If  $\tilde{f}= (j,j)$, then $H= (D_{4mr}^*/C_{2m}, D^*_{4nr}/C_{2n})_s$ (Family 11)
\end{itemize}

If $r$ is odd, we can suppose without loss of generality that $s$ is odd; otherwise, we can replace $s$ with $r-s$ obtaining  a conjugate group (see~\cite[page 803]{mecchia-seppi})

\begin{itemize}
\item Suppose  $\tilde{f}= (e^{i\pi/2mr}, e^{is\pi/2nr})$. We have both $e^{i\pi/2mr} \notin C_{2mr}$ and $e^{is\pi/2nr} \notin C_{2nr}$, so $L'=C_{4mr}, R'=C_{4nr}$, and we then find $L_K'= C_{2m} $, $R_K'= C_{2n} $ and $s'= s$, so $H=(C_{4mr}/C_{2m}, C_{4nr}/C_{2n})_s$ (Family 1).

\item Suppose $\tilde{f}= (e^{i\pi/2m}, e^{i(s+1)\pi/2n})$; since $s+1$ is even, we have $e^{i\pi/2m} \notin C_{2mr}$ and $e^{i(s+1)\pi/2n}= e^{r((s+1)/2)2i\pi/2nr} \in C_{2nr}$, and therefore $L'=C_{4mr}, R'=C_{2nr}$. We then find $L_K'= C_{4m} $, $R_K'= C_{2n}$. 
To find $s'$, we need to find an element of the form $(e^{2i\pi /4mr}, y)$ in $H$. We know that $a:=(e^{r 2i\pi/4mr}, e^{r((s+1)/2)2i\pi /2nr})$ and $b:=(e^{2(2i\pi) /4mr}, e^{s(2i\pi)/2nr})$ are in $H$, so, by writing $r=2k+1$ ($r$ is odd), we have $ab^{-k}=(e^{2i\pi /4mr},e^{ ((r(s+1)-2ks)/2) 2i\pi/4nr})= (e^{2i\pi /4mr},e^{((r+s)/2) 2i\pi  /2nr}) \in H$. We obtain $s'=(s+r)/2$, and  $H= (C_{4mr}/C_{4m}, C_{2nr}/C_{2n})_{(r+s)/2}$ (Family 1).

\item Suppose $\tilde{f}= (e^{i\pi (1+r)/2mr}, e^{i\pi (rs+r+s)/2nr})$; $r+1$ is even and $rs+r+s$ is odd, so $L'=C_{2mr}$, $R'=C_{4nr}$. We then find $L_K'= C_{2m} $, $R_K'=C_{4n} $ and $s'= 2s\% r$, so $H= (C_{2mr}/C_{2m}, C_{4nr}/C_{4n})_{2s\% r}$ (Family 1).

\item If  $\tilde{f}= (j,j)$, the computation is the same as in the  case of $r$ even, and  $H= (D_{4mr}^*/C_{2m}, D^*_{4nr}/C_{2n})_s$ (Family 11).

\end{itemize}

\begin{table}
\caption{Fibrations of the quotients: general case}
\label{fibrations-one}
\begin{tabular}{|l|l|l|l|c|}
	\hline
    group & extension (family) & base  & local & $e$\\
     &  & orbifold & invariants & \\
    \hline 
    \hline
    $(C_{2mr}/C_{2m}, C_{2nr}/C_{2n})_s$
     & 3 extensions in Family 1 & cf. Table~\ref{seifert-fam-1} &  cf. Table~\ref{seifert-fam-1} &  \\
 \small{$(*_1)$}
     &  depending on $m,n,r$ and $s$ & & &\\
     & (see Subsection~\ref{extensions}) & & &\\
     & $(D^*_{4mr}/C_{2m}, D^*_{4nr}/C_{2n})_s$ ($11$) & $D^2(;nr,nr)$ & cf. Table~\ref{seifert-fam-1} & -$\frac{m}{nr}$\\
    \hline
    $(C_{mr}/C_{m}, C_{nr}/C_{n})_s$
     & $(C_{2mr}/C_{m}, C_{2nr}/C_{n})_s$ ($1'$) &   $S^2(nr,nr)$ &  cf. Table~\ref{seifert-fam-1-prime} &  -$\frac{m}{nr}$\\
 \small{$(*_2)$}
     & $(C_{mr}/C_{2m}, C_{nr}/C_{2n})_{s\%(r/2)}$ ($1$) & $S^2(nr/2,nr/2)$   & cf. Table~\ref{seifert-fam-1} & -$\frac{4m}{nr}$ \\
     & $(C_{2mr}/C_{m}, C_{2nr}/C_{n})_{s+r}$ ($1'$) & $S^2(nr,nr)$ &  cf. Table~\ref{seifert-fam-1-prime} & -$\frac{m}{nr}$ \\
     & $(D^*_{2mr}/C_{m}, D^*_{2nr}/C_{n})_s$ ($11'$) & $D^2(;nr/2,nr/2)$ & cf. Table~\ref{seifert-fam-1-prime} & -$\frac{m}{nr}$\\
    \hline
    $(C_{2m}/C_{2m},D^*_{4n}/D^*_{4n})$ &
       $(C_{2m}/C_{2m}, D^*_{8n},D^*_{8n})$ (2) & $S^2(2,2,2n)$ & $ \frac{m}{2}, \frac{m}{2}, \frac{m}{2n}$ & -$\frac{m}{2n}$
       \\[0.5ex]
     \small{$\gcd(2n,m) = 1$}
      & $(C_{4m}/C_{4m},D^*_{4n}/D^*_{4n})$ (2) & $S^2(2,2,n)$ & $\frac{2m}{2}, \frac{2m}{2}, \frac{2m}{n}$ & -$\frac{2m}{n}$
      \\[0.5ex]
     & $(D_{4m}^*/D_{4m}^*,D^*_{4n}/D^*_{4n})$ (10) & \small{$n$ even:} $D^2(; 2, 2, n)$  & $\frac{m}{2}, \frac{m}{2}, \frac{m}{n}$ & -$\frac{m}{2n}$\\[0.5ex]

     & & \small{$n$ odd:}  $D^2(2;n)$  & $\frac{m}{2}, \frac{m}{n}$
     & -$\frac{m}{2n}$\\[0.5ex]
     & $(C_{4m}/C_{2m}, D_{8n}^*/D_{4n}^*)$ (4) & $S^2(2,2,2n)$ & $\frac{m}{2}, \frac{m+1}{2}, \frac{m+n}{2n}$
     & -$\frac{m}{2n}$ \\[0.5ex]
     & $(D^*_{4m}/C_{2m}, D_{8n}^*/D_{4n}^*)$ (13bis) & \small{$n$ odd:} $D^2(; 2, 2, n)$  &$\frac{m}{2}, \frac{m}{2}, \frac{m}{n}$ & -$\frac{m}{2n}$ \\[0.5ex]
     &  & \small{$n$ even:} $D^2(2;n) $ & $\frac{m}{2}, \frac{m}{n}$
     & -$\frac{m}{2n}$ \\
     \hline
    $(C_{4m}/C_{2m}, D_{4n}^*/C_{2n})$
     & $(C_{4m}/C_{2m}, D_{8n}^*/C_{4n})$ (3) & $S^2(2,2,2n) $ & $\frac{m+1}{2}, \frac{m+1}{2}, \frac{m}{2n}$
     & -$\frac{m}{2n}$ \\[0.5ex]
     \small{$\gcd(m,n) = 1$}
     & $(C_{4m}/C_{4m},D_ {4n}^*/D_ {4n}^*)$ (2) & $S^2(2,2,n)$ & $\frac{2m}{2}, \frac{2m}{2}, \frac{2m}{n}$
     & -$\frac{2m}{n}$\\[0.5ex]
    \small{and $m$ even} & $(D^*_{8m}/D^*_{4m}, D_{4n}^*/C_{2n})$ (13) & $D^2(2;n) $ & $\frac{m+1}{2}, \frac{m}{n}$
     & -$\frac{m}{2n}$ \\[0.5ex]
     & $(C_{4m}/C_{2m}, D_{8n}^*/D^*_{4n})$ (4) & $S^2(2,2,2n)$ & $ \frac{m}{2}, \frac{m+1}{2}, \frac{m+n}{2n}$
     & -$\frac{m}{2n}$ \\[0.5ex]
     & $(D^*_{8m}/C_{2m}, D_{8n}^*/C_{4n})$ (11) & $D^2(;2n,2n)$ &  cf. Table~\ref{seifert-fam-1} &  -$\frac{m}{2n}$\\
    \hline
    $(C_{2m}/C_{2m},T^*/T^*)$
     & $(C_{2m}/C_{2m},O^*/O^*)$ (7) & $S^2(2,3,4)$ & $\frac{m}{2}, \frac{m}{3}, \frac{m}{4}$ & -$\frac{m}{12}$ \\[0.5ex]
     \small{$\gcd(m,6) = 1$}
     & $(C_{4m}/C_{4m}, T^*/T^*)$ (5)  & $S^2(2,3,3)$ & $ \frac{2m}{2}, \frac{2m}{3}, \frac{2m}{3} $ & -$\frac{m}{3}$ \\[0.5ex]
     & $(D^*_{4m}/D^*_{4m},T^*/T^*)$ (14)  & $D^2(3;2)$  & $\frac{m}{2}, \frac{m}{3}$ & -$\frac{m}{12}$ \\[0.5ex]
     & $(C_{4m}/C_{2m}, O^*/T^*)$ (8) & $S^2(2,3,4)$ & $ \frac{m+1}{2}, \frac{m}{3}, \frac{m+2}{4}$ & -$\frac{m}{12}$\\[0.5ex]
     & $(D^*_{4m}/C_{2m}, O^*/T^*)$ (16)& $D^2(; 2, 3, 3)$ & $ \frac{m}{2}, \frac{m}{3}, \frac{m}{3}$ & -$\frac{m}{12}$ \\
    \hline
    $(C_{6m}/C_{2m}, T^*/D_8^*)$
     & $(C_{12m}/C_{4m}, T^*/D_8^*)$ (6) & $S^2(2,3,3)$ & $ \frac{2m}{2}, \frac{2m+1}{3}, \frac{2m+2}{3}$ & -$\frac{m}{3}$\\[0.5ex]
    \small{$m$ odd and $3|n$}
     & $(D^*_{12m}/C_{2m}, O^*/D_8^*$ (18) & $D^2(; 2, 3, 3) $ & $\frac{m}{2}, \frac{m+1}{3}, \frac{m+2}{3}$ & -$\frac{m}{12}$\\
    \hline
    $(C_{2m}/C_{2m},O^*/O^*)$
     & $(C_{4m}/C_{4m},O^*/O^*)$ (7) & $S^2(2,3,4)$ & $\frac{2m}{2}, \frac{2m}{3}, \frac{2m}{4}$ & -$\frac{m}{6}$\\[0.5ex]
     \small{$\gcd(m,6) = 1$}
     & $(D_{4m}/D_{4m}, O^*/O^*)$ (15) & $D^2(; 2, 3, 4)$ & $ \frac{m}{2}, \frac{m}{3}, \frac{m}{4}$ & -$\frac{m}{24}$\\
    \hline
    $(C_{2m}/C_{2m},I^*/I^*)$
     & $(C_{4m}/C_{4m},I^*/I^*)$ (9) & $S^2(2,3,5)$ & $\frac{2m}{2}, \frac{2m}{3}, \frac{2m}{5}$ & -$\frac{m}{15}$\\[0.5ex]
            \small{$\gcd(m,30) = 1$}
     & $(D^*_{4m}/D^*_{4m}, I^*/ I^*)$ (19) & $D^2(; 2, 3, 5) $ & $\frac{m}{2}, \frac{m}{3}, \frac{m}{5}$ & -$\frac{m}{60}$\\
    \hline
     \multicolumn{5}{|l|}{\small{$(*_1)$: if $\gcd(r,s)=1,$ $\gcd(m,n)=1$, $mn$ even or $r$ odd,}} \\
     \multicolumn{5}{|l|}{and $\gcd(n-sm,2mnr)= \gcd(n+sm,2mnr)$ } \\
     \hline
     \multicolumn{5}{|l|}{\small{$(*_2)$: if $\gcd(r,s)=1,$ $r$ even, $m$ odd, $n$ odd }}\\
     \multicolumn{5}{|l|}{\small{$\gcd(m,n)=1$ and $\gcd(n-sm,mnr)= \gcd(n+sm,mnr)$ }} \\
     \hline
\end{tabular}
\end{table}


\begin{table}
\caption{Fibration of the quotients: small indices (first part)}
\label{fibrations-two}
\begin{tabular}{|l|l|l|l|c|}
	\hline
     group  & extension (family) & base  & local & $e$\\
     &  & orbifold & invariants & \\
    \hline  
    \hline
    $(C_{2m}/C_{2m}, C_{2n}/C_{2n})$
     & $(C_{4m}/C_{4m},C_{2n}/C_{2n})$  (1) & $S^2(n,n)$ & cf. Table~\ref{seifert-fam-1} & -$\frac{m}{n}$  \\
 \small{$\gcd(n,m)=1$}
     & $(D^*_{4m}/D^*_{4m}, C_{2n}/ C_{2n})$ (2bis) & \small{$n$ even:} $D^2(n)$ & $\frac{m}{n}$& -$\frac{m}{n}$ \\
     & &\small{$n$ odd}: $\mathbb{R}P^2(n)$ & $\frac{m}{n}$& -$\frac{m}{n}$  \\
     & $(C_{4m}/C_{2m}, C_{4n}/C_{2n})$ (1) & $S^2(n,n)$ & cf. Table~\ref{seifert-fam-1} & -$\frac{m}{n}$ \\
     & $(D^*_{4m}/C_{2m}, D^*_{4n}/C_{2n})$ (11) & $D^2(;n,n)$ &  cf. Table~\ref{seifert-fam-1} & -$\frac{m}{n}$ \\
     & $(C_{4m}/C_{2m}, D^*_{4n}/C_{2n})$ (3) & $S^2(2,2,n)$ & $\frac{m+1}{2}, \frac{m+1}{2}, \frac{m}{n}$& -$\frac{m}{n}$\\
     & $(C_{2m}/C_{2m},C_{4n}/C_{4n})$ (1) &  $S^2(2n,2n)$& cf. Table~\ref{seifert-fam-1}  &-$\frac{m}{n}$\\
     & $(C_{2m}/C_{2m},D^*_{4n}/D^*_{4n})$ (2) & $S^2(2,2,n)$ & $\frac{m}{2}, \frac{1}{2}, \frac{m}{n}$   &-$\frac{m}{n}$ \\
     & $(D^*_{4m}/C_{2m}, C_{4n}/C_{2n})$ (3bis) & \small{$n$ odd} $D^2(n)$ & $\frac{m}{n}$ &-$\frac{m}{n}$\\
     & & \small{$n$ even}  $\mathbb{R}P^2(n)$ & $\frac{m}{n}$&-$\frac{m}{n}$\\
    \hline
    $(C_{4m}/C_{2m}, C_{4n}/C_{2n})$ 
     & $(C_{4m}/C_{4m},C_{4n}/C_{4n})$ (1) & $S^2(2n,2n)$ & cf. Table~\ref{seifert-fam-1} & -$\frac{2m}{n}$\\
 \small{$\gcd(m,n)=1$ and}
     & $(D^*_{8m}/D^*_{4m}, C_{4n}/C_{2n})$ (4bis) & $D^2(2n)$ & $\frac{m+n}{2n}$& -$\frac{m}{2n}$\\
      \small{$mn$ even} & $(D_{8m}^*/C_{2m}, D_{8n}^*/C_{2n})$ (11) & $D^2(;2n,2n)$  & cf. Table~\ref{seifert-fam-1}&-$\frac{m}{2n}$\\
     & $(C_{8m}/C_{2m}, C_{8n}/C_{2n})$ (1) &  $S^2(4n,4n)$ & cf. Table~\ref{seifert-fam-1} & -$\frac{m}{2n}$\\
     & $(C_{4m}/C_{2m}, D^*_{8n}/D^*_{4n})$ (4) & $S^2(2,2,2n)$  & $\frac{m}{2}, \frac{m+1}{2}, \frac{m+n}{2n}$ & -$\frac{m}{2n}$\\
    \hline
    $(C_{2}/C_{2},C_{2n}/C_{2n})$ 
     & $(C_{4}/C_{4},C_{2n}/C_{2n})$ (1) & $S^2(n,n)$ & cf. Table~\ref{seifert-fam-1} & -$\frac{4}{n}$\\
     & $(C_{2}/C_{2},C_{4n}/C_{4n})$ (1) & $S^2(2n,2n)$ & cf. Table~\ref{seifert-fam-1} & -$\frac{1}{n}$\\
     & $(C_{2}/C_{2}, D^*_{4n}/D^*_{4n})$ (2) & $S^2(2,2,n)$ & $ \frac{1}{2}, \frac{1}{2}, \frac{1}{n}$ & -$\frac{1}{n}$\\
     & $(C_{4}/C_{2}, C_{4n}/C_{2n})$ (1) & $S^2(2n,2n)$ & cf. Table~\ref{seifert-fam-1} & -$\frac{1}{n}$\\
     & $(C_{4}/C_{2}, D^*_{4n}/C_{2n})$ (3) & $S^2(2,2,n)$  & $\frac{2}{2}, \frac{2}{2}, \frac{1}{n}$& -$\frac{1}{n}$\\
    \hline
      $(C_{2}/C_{2},C_{2}/C_{2})$ 
     & $(C_{4}/C_{2}, C_{4}/C_{2})$ (1)  & $S^2(2,2)$ & $\frac{4}{2}, -\frac{2}{2}$& -1\\
     & $(C_{4}/C_{4},C_{2}/C_{2})$ (1) & $S^2$ & cf. Table~\ref{seifert-fam-1}& -4\\
     & $(C_{2}/C_{2},C_{4}/C_{4})$ (1) & $S^2(2,2)$ & cf. Table~\ref{seifert-fam-1} & -1 \\
     \hline
    $(C_{2m}/C_{m}, C_{2n}/C_{n})$ 
     & $(C_{2m}/C_{2m},C_{2n}/C_{2n})$ (1) & $S^2(n,n)$ & cf. Table~\ref{seifert-fam-1} & -$\frac{2m}{n}$\\
    \small{$mn$ odd and }
     & $(C_{4m}/C_{m}, C_{4n}/C_{n})$ ($1'$) & $S^2(2n,2n)$  & cf. Table~\ref{seifert-fam-1-prime} & -$\frac{m}{2n}$\\
    \small{$\gcd(m,n)=1$}  & $(D^*_{4m}/C_{m}, D^*_{4n}/C_{n})$ ($11'$) &  $D^2(;n,n))$&  cf. Table~\ref{seifert-fam-1-prime} & -$\frac{m}{2n}$\\
     & $(C_{4m}/C_{m}, D^*_{4n}/C_{n})$ (34) & $S^2(2,2,n)$ & $\frac{m}{2}, \frac{m+1}{2}, \frac{m+n}{n}$  & -$\frac{m}{2n}$\\
     & $(D^*_{4m}/C_{m}, C_{4n}/C_{n})$ (34bis) & $D^2(n;)$  & $\frac{(m+n)/2}{n}$ &-$\frac{m}{2n}$\\   
    \hline
    $(C_{2}/C_1, C_{2n}/C_{n})$ 
     & $(C_{2}/C_{2},C_{2n}/C_{2n})$ (1) & $S^2(n,n)$ & cf. Table~\ref{seifert-fam-1}  & -$\frac{2}{n}$ \\
     \small{$n$ odd}  & $(C_{4}/C_{1}, C_{4n}/C_{n})$ ($1'$) &  $S^2(2n,2n)$ &  cf. Table~\ref{seifert-fam-1-prime} & -$\frac{1}{2n}$ \\
     & $(C_{4}/C_{1}, D^*_{4n}/C_{n})$ (34) & $S^2(2,2,n)$  & $\frac{1}{2}, \frac{2}{2}, \frac{(n+1)/2}{n}$& -$\frac{1}{2n}$\\
    \hline
    $(C_{2}/C_{1}, C_{2}/C_{1})$ 
     & $(C_{4}/C_{1}, C_{4}/C_{1})_1$ ($1'$) & $S^2(2,2)$ & $\frac{2}{2}, \frac{1}{2}$  & -$\frac{1}{2}$ \\
     & $(C_{2}/C_{2},C_{2}/C_{2})$ (1) & $S^2$ &  & -2\\
    \hline
\end{tabular}
\end{table}


\begin{table}
\caption{Fibration of the quotients: small indices (second part)}
\label{fibrations-three}
\begin{tabular}{|l|l|l|l|c|}
	\hline
     group (family) & extension  & base  & local & $e$ \\
     &  & orbifold & invariants & \\
    \hline
    \hline
   $(C_{2}/C_{2},D_{4n}^*/D_{4n}^*)$  
     & $(C_{4}/C_{4}, D_{4n}^*/D_{4n}^*)$ (2) & $S^2(2,2,n)$ & $\frac{2}{2}, \frac{2}{2}, \frac{2}{n}$ & -$\frac{2}{n}$\\[0.5ex]
     & $(C_{2}/C_{2}, D_{8n^*}/D_{8n^*})$  (2) & $S^2(2,2,2n)$ & $\frac{1}{2}, \frac{1}{2}, \frac{1}{2n}$ & -$\frac{1}{2n}$\\[0.5ex]
     & $(C_{4}/C_{2}, D^*_{8n}/D^*_{4n})$ (4) & $S^2(2,2,2n)$ & $ \frac{1}{2}, \frac{2}{2}, \frac{n+1}{2n}$ &-$\frac{1}{2n}$\\[0.5ex]
    \hline
    $(C_{2m}/C_{2m}, D_{8}/D_{8})$  
     & $(C_{4m}/C_{4m},D_{8}^*/D_{8}^*)$ (2) & $S^2(2,2,2)$   & $\frac{2m}{2}, \frac{2m}{2}, \frac{2m}{n}$& -$m$ \\[0.5ex]
     \small{$m$ odd}
     & $(D^*_{4m}/D^*_{4m},D_{8}^*/D_{8}^*)$ (10) & $D^2(;2,2,2)$ & $\frac{m}{2}, \frac{m}{2}, \frac{m}{2}$ & -$\frac{m}{4}$\\[0.5ex]
     & $(C_{2m}/C_{2m},D_{16}^*/D_{16}^*)$ (2) & $S^2(2,2,4)$ & $\frac{m}{2}, \frac{m}{2}, \frac{m}{4}$ &-$\frac{m}{4}$\\[0.5ex]
     & $(C_{4m}/C_{2m}, D^*_{16}/D^*_{8})$ (4) & $S^2(2,2,n)$ & $\frac{m}{2}, \frac{m+1}{2}, \frac{m+2}{4}$&-$\frac{m}{4}$\\[0.5ex]
     & $(D^*_{4m}/C_{2m}, D^*_{16}/D^*_{8})$ (13bis) & $D^2(2;2)$ & $\frac{m}{2}, \frac{m}{2}, \frac{m}{2}$ &-$\frac{m}{4}$\\[0.5ex]
    \hline
   $(C_{2}/C_{2},D_{8}/D_{8})$ 
     & $(C_{4}/C_{4}, D_{8}^*/D_{8}^*)$  (2) & $S^2(2,2,2)$ & $\frac{2}{2}, \frac{2}{2}, \frac{2}{2}$ & -1\\[0.5ex]
     &  $(C_{2}/C_{2},D_{16}^*/D_{16}^*)$ (2) & $S^2(2,2,4)$  & $\frac{1}{2}, \frac{1}{2}, \frac{1}{4}$& -$\frac{1}{4}$\\[0.5ex]
     & $(C_{4}/C_{2}, D^*_{16}/D^*_{8})$ (4)  & $S^2(2,2,4)$ & $\frac{1}{2}, \frac{2}{2}, \frac{3}{4}$  & -$\frac{1}{4}$\\[0.5ex]
    \hline
     $(C_2/C_2,T^*/T^*)$ 
     &  $(C_4/C_4,T^*/T^*)$(5) & $S^2(2,3,3)$ & $\frac{2}{2}, \frac{2}{3}, \frac{2}{3}$ & -$\frac{1}{3}$\\[0.5ex]
     &$(C_2/C_2,O^*,O^*)$ (7) & $S^2(2,3,4)$  & $\frac{1}{2}, \frac{1}{3}, \frac{1}{4}$ & -$\frac{1}{12}$ \\[0.5ex]
     & $(C_4/C_2, O^*/T^*)$, (8) & $S^2(2,3,4)$ & $\frac{2}{2}, \frac{1}{3}, \frac{3}{4}$ & -$\frac{1}{12}$ \\[0.5ex]
    \hline
    $(C_2/C_2, O^*/O^*)$  & $(C_4/C_4,O^*/O^*)$ (7) & $S^2(2,3,4)$  & $\frac{2}{2}, \frac{2}{3}, \frac{2}{4}$& -$\frac{1}{6}$\\[0.5ex]
    \hline
     $(C_2/C_2,I^*/I^*)$  & $(C_4/C_4,I^*/I^*)$ (9) & $S^2(2,3,5)$  & $\frac{2}{2}, \frac{2}{3}, \frac{2}{5}$ & -$\frac{1}{15}$\\[0.5ex]
    \hline
\end{tabular}
\end{table}

\subsection{Seifert fibrations of the quotient orbifolds.} \label{geometric}

For any involution of a spherical 3-manifold, the group $\Phi(H)$  leaves invariant the Hopf fibration of $S^3$, thus the Hopf fibration induces a Seifert fibration on the quotient orbifold $(S^3/G)/\langle f \rangle=S^3/\Phi(H).$  We can use the formulae given in \cite{mecchia-seppi} to compute the invariants of the Seifert fibration of $(S^3/G)/\langle f \rangle $.  We use the same notation used in \cite{mecchia-seppi}. For all groups but those in Families 1,$1'$, 11 and $11'$, the local invariants are given by closed formulae, therefore for these groups, the values of the local invariants are directly inserted in Tables~\ref{fibrations-one}, \ref{fibrations-two} and \ref{fibrations-three}. For Families 1,$1'$, 11 and $11'$ the computation is more complicated and, in many cases,  the results are difficult to insert in a single cell of a table, so we report in Table~\ref{seifert-fam-1} and \ref{seifert-fam-1-prime} the details of the formulae that give the invariants for these four families. The last column of Tables~\ref{fibrations-one}, \ref{fibrations-two} and \ref{fibrations-three} contains the  Euler classes of the fibrations.

We can deduce many geometric features of the involution from the invariants of the fibration. For the sake of simplicity, we denote $(S^3/G)/\langle f \rangle=S^3/\Phi(H)$ by $\OO_f$ and its base orbifold  by $\BB_f.$ 
\bigskip

\subsection{Hyperelliptic involutions.} \label{hyperelliptic}

If the underlying topological space of $\OO_f$ is $S^3$, then the underlying topological space of the base orbifold $\BB_f$ is either  $D^2$ or $S^2$, see  \cite[Section 5]{dunbar}. We recall that, if $m/n$ is a local invariant of a cone point of $\BB_f$, the $\gcd(m,n)$ is the index of singularity of the fibre. If we replace $m/n$ with  $(m/\gcd(m,n))/(n/\gcd(m,n))$ the underlying topological space of $\OO_f$ does not change. We call $(m/\gcd(m,n))/(n/\gcd(m,n))$ the \textit{reduced local invariant.}

If the underlying topological space of $\BB_f$ is a disk, the base orbifold is denoted by $$D^2(a_1,a_2,\dots;b_1,b_2,\dots)$$

where the integers appearing before the semicolon represent the cone points in the interior of the disk.

The  Seifert orbifolds appearing as the quotient of involutions in Tables~\ref{fibrations-one}, \ref{fibrations-two} and \ref{fibrations-three} admit at most one cone point in the interior of the disk. 

If  $\BB_f$ has no cone points, the underlying topological space of $\OO_f$ is $S^3,$ see \cite[Proposition 2.10]{dunbar2}.

If we have one cone point, the underlying topological space of $\OO_f$ is $S^3$ if and only if the reduced local invariant of the cone point is an integer, see \cite[Proposition 2.11]{dunbar2}. This situation corresponds to a fibre that is singular but not exceptional.

We consider now the case when $\BB_f$ is a sphere.  The underlying topological space of $\OO_f$ is $S^3$ if and only if the fibration obtained replacing the local invariants with the reduced ones gives a Seifert fibration of the sphere.

The  Seifert fibrations for $S^3$ are well known (see \cite{seifert} or also \cite[Sections 2.5 and 4.1]{MecchiaSeppi3}): they have 
 base orbifold $S^2(u,v)$ for $u,v\geq 1$ two coprime integers, hence at most two exceptional fibres.  The local invariants are given by (the classes modulo 1 of) $\bar v/u$ and $\bar u/v$ where $u\bar u+v \bar v=1,$ and the Euler class is $\pm 1/uv$. If $u$ or $v$ equals $1$, we mean that the corresponding point in the base orbifold is regular and hence there is no local invariant to associate (the above formula would indeed give 0 as output). In particular, for $u=v=1$ we obtain the \textit{Hopf fibration}.

If the group $G$ is in Family 1 or $1'$, then the invariants of the fibration are difficult to compute, hence we use another method to determine whether the underlying space is $S^3$. Table~\ref{seifert-fam-1} and \ref{seifert-fam-1-prime} describes the underlying topological space  of $S^3/G$ as a lens space $L(l, *)$, where $l$ is an integer depending on $m,n,r$ and $s$. Therefore, the underlying space is $S^3$ if and only if  $l=1$. The integer $l$ can always be determined provided that it is analysed on a case-by-case basis, but a global analysis might be difficult. Here we prove that the extensions in Family 1 and $1'$ are not associated with hyperelliptic involutions except $(C_4/C_1; C_4/C_1).$ The following remark is useful during the analysis.

\begin{remark} \label{computation-f}

If $\tilde{G}=(C_{2mr}/C_{2m}, C_{2nr}/C_{2n})_s$ is a group of Family 1 such that the underlying space of the quotient $S^3/G$ is $S^3$ we have $l=1$ in Table~\ref{seifert-fam-1}; this implies  $2m'n'r=l_1l_2b_1b_2$, with $h=\gcd(m,n),$ $ m'=m/h$ and $n'=n/h$. Let $p\neq 2$ be a prime number dividing $m'$. Then $p|b_1b_2$, so $p|b_1$ or $p|b_2$, and in both cases $p|n'$, that is impossible. Therefore $m'=2^k$; if we suppose that $k\geq 1$  then $l_1=l_2=1$ e $2|b_1b_2$, so $2|n'$, which is impossible. We obtain that if the underlying space of the quotient $S^3/G$ is $S^3$ then necessarily $m'=1$.

Now we will show that the same holds for $n'.$ We recall we represent $S^3$ as the set of unit quaternions. Let $\gamma:S^3 \rightarrow S^3 $ be the orientation-reversing isometry of $S^3$ such that $\gamma(h)=\bar{h}$.  We have that 
$\gamma^{-1}\tilde{G}\gamma=(C_{2nr}/C_{2n}, C_{2mr}/C_{2m})_{\hat{s}}$; the parameter $\hat{s}$ might be different from $s$ (see \cite[Section 3.2]{mecchia-seppi2}) but for our purpose is not relevant. The quotient orbifold $S^3/\Phi(\gamma^{-1}\tilde{G}\gamma)$ is diffeomorphic by an orientation-preserving diffeomorphism to $S^3/G$ and we obtain again $S^3$ as underlying topological space. This implies $n'=1.$

The same result can be obtained for the Family $1'$.
    
\end{remark}

For example, we consider the Family $1'$ extensions in the general case, that is  $\tilde{G}=(C_{mr}/C_{m}, C_{nr}/C_{n})_s$ in Table~\ref{fibrations-one}.  We have four involutions up to conjugacy; the fourth one in the table gives as quotient a Seifert fibred orbifold whose base orbifold is a disk without cone points so the underlying topological space is the 3-sphere and the involution is hyperelliptic. The other involutions are given by extensions of Family 1 and $1'$ and a more careful analysis is needed.

Consider $f$ the involution given by the extension $(C_{2mr}/C_{m}, C_{2nr}/C_{n})_s$ of Family $1'$ and suppose that the underlying topological space of $\OO_f$  is $S^3$. We recall that since $S^3/G$ is a manifold we have $\gcd(n,m)=1$ thus $m=m'$ and $n=n'$; by Remark~\ref{computation-f} we obtain $m=n=1.$ Since  $\gcd(r,s)=1$ and $r$ even then $\gcd(1-s,1+s, r)=2$.  Since $S^3/G$  is a manifold,  the singular set is trivial and this implies that  $\gcd((1-s)/2,r/2)=1$ and $\gcd((1+s)/2,r/2)=1,$ see Table~\ref{seifert-fam-1-prime}. One between $(1-s)/2$ and $(1+s)/2$ is even so  we obtain that $r/2$ is odd.

  The underlying topological space  of $\OO_f$ is a lens space $L(l, *)$ where 

$$l=\frac{2r}{2\gcd((1-s)/2,r)\gcd((1-s)/2,r)}.$$

Considering that  $\gcd((1-s)/2,r/2)=1=\gcd((1+s)/2,r/2)=1$ and $r/2$ odd, we obtain that $l=(r/2)>1$ (in the general case $r>2$). We can conclude that $f$ cannot be hyperelliptic. 

If the involution is given by the extension $(C_{2mr}/C_{m}, C_{2nr}/C_{n})_{s+r}$ the analysis is very similar and the involution cannot be hyperelliptic.

Finally,  we consider the involution given by $(C_{mr}/C_{2m}, C_{nr}/C_{2n})_{s\%(r/2)}$ in Family 1. To use formulae in Table~\ref{seifert-fam-1} we have to see the group in the following way: $$(C_{2m(r/2)}/C_{2m}, C_{2n(r/2)}/C_{2n})_{\hat{s}}$$
where $\hat{s}$ is odd (if $s\%(r/2)$ is even we can replace it with $(r/2)-s\%(r/2)).$

Again we obtain $m=n=1$, $\gcd((1-s)/2,r/2)=1=\gcd((1+s)/2,r/2)=1$ and $r/2$ odd.  If the underlying topological space of $\OO_f$ is the lens space $L(l, *)$ we have that $$l=\frac{2(r/2)}{l_1l_2b_1b_2}.$$ 

Since $r/2$ is odd $l_1=l_2=1$. Since $\hat{s}$ is odd we have that $\gcd(1-\hat{s},1+\hat{s}, 2(r/2))=2$ and $b_1=\gcd((1-\hat{s})/2, r/2)=1$ and $b_2=\gcd((1+\hat{s})/2, r/2)=1$. We obtain that $l=(r/2)>1.$

For the other groups, similar computations can be carried out. The results are collected in the last column of Tables~\ref{general quotient} and \ref{small index quotient}. If ``Y'' appears the involution is hyperelliptic, while ``N'' means that it is not; if a condition occurs in the column, the involution is hyperelliptic if and only if the condition is verified. We note that the group $(C_2/C_1, C_2/C_1)$ is the trivial one and the corresponding quotient is the 3-sphere. It is well known that in this case, we have two conjugacy classes of involutions: one containing involutions acting freely and the other one consisting of hyperelliptic involutions.

\begin{table}
\caption{Seifert fibrations for Families $1$ and $11$.}
\label{seifert-fam-1}
\begin{tabular}{|l|l|}
\hline
\multicolumn{2}{|c|}{} \\
\multicolumn{2}{|c|}{$G=\Phi((C_{2mr}/C_{2m},C_{2nr}/C_{2n})_s)$ (Family $1$)}\\
\multicolumn{2}{|c|}{and}\\
\multicolumn{2}{|c|} {$G=\Phi((D^*_{4mr}/C_{2m},D^*_{4nr}/C_{2n})_s)$ (Family $11$)} \\
\multicolumn{2}{|c|}{} \\
\hline
\multicolumn{2}{|c|}{} \\
\multicolumn{2}{|l|}{We define: $h=\gcd(m,n),$ $m'=\frac{m}{h},$ $n'=\frac{n}{h},$ $a=\gcd(n'-sm',n'+sm',2m'n'r),$}\\
\multicolumn{2}{|l|}{  $b_1=\gcd(\frac{n'-sm'}{a},\frac{2m'n'r}{a}),$ $b_2=\gcd(\frac{n'+sm'}{a},\frac{2m'n'r}{a}).$} \\  
\multicolumn{2}{|c|}{} \\

\multicolumn{2}{|l|}{Remark: W.L.O.G. we assume $s$ odd}\\
\multicolumn{2}{|c|}{} \\
 \hline
&\\
if $n'm'$ is even,  we define: & in both cases we define:  \\
$\nu$ minimal positive integer s.t. $\gcd(\frac{n'}{a\nu },a)=1$ &  $d=\frac{\nu^2 a(n'+sm')+2n'm'r}{f_2a\nu b_2}$  \\ 
 $l_1=l_2=1$ & $g=\frac{\nu^2 a(n'-sm')-2n'm'r}{l_1a\nu b_1}$\\
& $l=\frac{2m'n'r}{l_1l_2b_1b_2}$ \\
if $n'm'$ is odd,  we define: &   $\overline{g}$ s.t. $g\overline{g}\equiv 1 \, \mod l$\\
$\nu$  minimal positive integer s.t.  $\gcd(\frac{2n'}{a\nu },\frac{a}{2})=1$  & $\overline{c}$ s.t. $\left(\nu s+r\frac{2n'}{a\nu }\right)\overline{c}\equiv 1 \, \mod n'r$ \\
																							$l_i=\left\{\begin{array}{ll}
																							2 &  \text{ if } \frac{r}{b_i} \text{ is even }\\
																							1 &  \text{ if } \frac{r}{b_i} \text{ is odd }
																							\end{array}\right.$ &   \\
&  \\
\hline
\multicolumn{2}{|c|}{} \\
\multicolumn{2}{|l|} {The orbifold $S^3/\Phi((C_{2mr}/C_{2m},C_{2nr}/C_{2n})_s)$ fibers over $S^2(nr,nr)$}\\

\multicolumn{2}{|l|} { with local invariants $\frac{d \overline{c}l_2b_2h}{nr}$ and $-\frac{g \overline{c}l_1b_1h}{nr}$ and Euler class $-\frac{2m}{nr}$.}\\
\multicolumn{2}{|l|} {The underlying topological space of $S^3/\Phi((C_{2mr}/C_{2m},C_{2nr}/C_{2n})_s)$ } \\
\multicolumn{2}{|l|} {is the lens space  $ L(l,d\overline{g})$.}\\
\multicolumn{2}{|l|} {The singular set  of $S^3/\Phi((C_{2mr}/C_{2m},C_{2nr}/C_{2n})_s)$ is a link with  } \\ 
\multicolumn{2}{|l|} {at most two components of singular index $l_2b_2 h$ and $l_1b_1 h$}\\
\multicolumn{2}{|l|} {(if the singular index is 1 the corresponding component } \\
\multicolumn{2}{|l|} {consists of non-singular points).}\\
\multicolumn{2}{|c|}{} \\
\hline
\multicolumn{2}{|c|}{} \\
\multicolumn{2}{|l|} {The orbifold $S^3/\Phi((D^*_{4mr}/C_{2m},D^*_{4nr}/C_{2n})_s)$ fibers over $D^2(;nr,nr)$}\\ \multicolumn{2}{|l|} { with local invariants   $\frac{d \overline{c}l_2b_2h}{nr}$ and $-\frac{g \overline{c}l_1b_1h}{nr}$ and Euler class $-\frac{m}{nr}$.}\\
\multicolumn{2}{|l|} {The underlying topological space of $S^3/\Phi((D^*_{4mr}/C_{2m},D^*_{4nr}/C_{2n})_s)$ }\\
\multicolumn{2}{|l|} { is the 3-sphere.}\\
\multicolumn{2}{|c|}{} \\
\hline
\end{tabular}
\end{table}

\begin{table}
\caption{Seifert fibrations for Families $1'$ and $11'$.}
\label{seifert-fam-1-prime}
\begin{tabular}{|ll|}
\hline
\multicolumn{2}{|c|}{} \\
\multicolumn{2}{|c|}{$G=\Phi((C_{mr}/C_{m},C_{nr}/C_{n})_s)$ (Family $1'$)}\\
\multicolumn{2}{|c|}{and}\\
\multicolumn{2}{|c|}{ $G=\Phi((D^*_{2mr}/C_{m},D^*_{2nr}/C_{n})_s)$ (Family $11'$)} \\
\multicolumn{2}{|c|}{} \\
\hline
&\\
We define:  & $h=\gcd(m,n)$ \\
						&  $m'=\frac{m}{h}$ \\
					 &  $n'=\frac{n}{h}$ \\ 
					 & $a=\gcd(n'-sm',n'+sm',m'n'r)$ \\
 & $b_1=\gcd(\frac{n'-sm'}{a},\frac{m'n'r}{a})$ \\
 &  $b_2=\gcd(\frac{n'+sm'}{a},\frac{m'n'r}{a})$ \\   
 &   $\nu $ minimal positive integer s.t.  $\gcd(\frac{2n'}{a\nu },\frac{a}{2})=1$ \\
& $d=\frac{\nu^2 a(n'+sm')+2n'm'r}{2a\nu b_2}$ \\ 
 & $g=\frac{\nu^2 a(n'-sm')-2n'm'r}{2a\nu b_1}$  \\
 & $l=\frac{m'n'r}{2b_1b_2}$ \\
 & $\overline{g}$ s.t. $g\overline{g}\equiv 1 \, \mod l$ \\
& $\overline{c}$ s.t. $\left(\nu s+r\frac{2n'}{a\nu }\right)\overline{c}\equiv 1 \, \mod n'r$ \\
&\\
\hline
\multicolumn{2}{|c|}{} \\
\multicolumn{2}{|p{14cm}|} {The orbifold $S^3/\Phi((C_{mr}/C_{m},C_{nr}/C_{n})_s)$ fibers over $S^2\left(\frac{nr}{2},\frac{nr}{2}\right)$ with local invariants $\frac{d \overline{c}b_2 h}{\frac{nr}{2}}$ and $-\frac{g \overline{c}b_1h}{\frac{nr}{2}}$ and Euler class $-\frac{2m}{nr}$.}\\
&\\
\multicolumn{2}{|p{14cm}|} {The underlying topological space of $S^3/\Phi((C_{mr}/C_{m},C_{nr}/C_{n})_s)$ is the lens space  $ L(l,d\overline{g})$.}\\
&\\
\multicolumn{2}{|p{14cm}|} {The singular set of $S^3/\Phi((C_{mr}/C_{m},C_{nr}/C_{n})_s)$ is a link with at most two components of singular index $b_2h$ and $b_1h$ (if the singular index is 1 the corresponding component consists of non-singular points).}\\
\multicolumn{2}{|c|}{} \\
\hline
\multicolumn{2}{|c|}{} \\
\multicolumn{2}{|p{14cm}|} {The orbifold $S^3/\Phi((D^*_{2mr}/C_{m},D^*_{2nr}/C_{n})_s)$ fibers over $D^2\left(;\frac{nr}{2},\frac{nr}{2}\right)$ with local invariants $\frac{d \overline{c}b_2 h}{\frac{nr}{2}}$ and  $-\frac{g \overline{c}b_1h}{\frac{nr}{2}}$ and Euler class $-\frac{m}{nr}$.}\\
&\\
\multicolumn{2}{|p{14cm}|} {The underlying topological space of $S^3/\Phi((D^*_{2mr}/C_{m},D^*_{2nr}/C_{n})_s)$ is the 3-sphere.}\\
\multicolumn{2}{|c|}{} \\
\hline
\end{tabular}
\end{table}
\bigskip

\subsection{Involutions acting freely.} 

If the Seifert fibration invariants of $\OO_f$ are explicit, it is easy to compute the number of components and to check whether $f$ acts freely (in this case $\OO_f$ is a manifold). 

We recall the following facts about Seifert fibrations:

\begin{itemize}
    \item If the underlying topological space of  $\BB_f$ is a disk, we have at least one singular fibre and $f$ does not act freely
    \item If the underlying topological space of  $\BB_f$ is a 2-sphere or the projective plane, the possible singular points of $\OO_f$ are contained in the fibres corresponding to the cone points of $\BB_f$; in particular, if $p/q$ is the local invariant of a cone point, the points of the fibre have singularity index $\gcd(p,q).$ If  $\gcd(p,q)=1$ the fibre consist of non-singular points. 
\end{itemize}

If the extension corresponding to the involution is not in Family 1 or $1'$, these two facts allow us to understand whether $f$ acts freely or, equivalently, when $\OO_f$ is a manifold.

For Family 1 and Family $1'$, we can use again the formulae in Table~\ref{seifert-fam-1} and \ref{seifert-fam-1-prime} to obtain this information. Indeed,  the situation seems too complicated to obtain a closed formula that fully describes the situation. 

For the small indices groups of Table~\ref{small index quotient}, the situation is clear also for groups in Family 1 and $1'$ so the table is complete.

\begin{table}
\caption{Geometric properties of involutions, general case}
\label{general quotient}
\begin{tabular}{|l|l|c|c|}
	\hline
    group (family) & extension (family) & free & hyperelliptic \\
 & & action &   \\
    \hline 
    \hline
    $(C_{2mr}/C_{2m}, C_{2nr}/C_{2n})$ (1)
    & 3 extensions in Family 1 &  &   \\
       \small{$(*_1)$} & depending on $m,n,r$ and $s$ &  & N \\
     & (see Subsection~\ref{extensions}) &  &  \\[0.5ex]\cline{2-4}

    & $(D^*_{4mr}/C_{2m}, D^*_{4nr}/C_{2n})_s$ (11)  & N & Y \\[0.5ex]
    \hline
    $(C_{mr}/C_{m}, C_{nr}/C_{n})_s$ ($1'$)
     & $(C_{2mr}/C_{m}, C_{2nr}/C_{n})_s$ ($1'$) &  & N \\
    \small{$(*_2)$}
     & $(C_{mr}/C_{2m}, C_{nr}/C_{2n})_{s\%(r/2)}$ (1) &  & N \\
     & $(C_{2mr}/C_{m}, C_{2nr}/C_{n})_{s+r}$ &  & N \\
     & $(D^*_{2mr}/C_{m}, D^*_{2nr}/C_{n})_s$ ($11'$) & N & Y \\
     \hline
    $(C_{2m}/C_{2m}, D^*_{4n}/D^*_{4n})$ (2) &
       $(C_{2m}/C_{2m}, D^*_{8n}/D^*_{8n})$ (2) & Y & N  
       \\
     \small{$\gcd(2n,m) = 1$}
      & $(C_{4m}/C_{4m}, D^*_{4n}/D^*_{4n})$ (2) & N & N
      \\
     & $(D^*_{4m}/D^*_{4m},  D^*_{4n}/D^*_{4n})$ (10) & N & \small{$n$ even}
     \\
     & $(C_{4m}/C_{2m}, D_{8n}^*/D_{4n}^*$) (4) & N & N
     \\
     & $(D^*_{4m}/C_{2m}, D_{8n}^*/D_{4n}^*)$ (13bis) & N & \small{$n$ odd}
     \\
     \hline
    $(C_{4m}/C_{2m}, D_{4n}^*/C_{2n})$ (3)
     & $(C_{4m}/C_{2m}, D_{8n}^*/C_{4n})$ (3) & N & N
     \\
     \small{$\gcd(m,n) = 1$ and $m$ even}
     & $(C_{4m}/C_{4m}, D_ {4n}^*/D_ {4n}^*)$ (2) & N & N
     \\
     & $(D^*_{8m}/D^*_{4m}, D_{4n}^*/C_{2n})$ (13)  & N & N
     \\
     & $(C_{4m}/C_{2m}, D_{8n}^*/D^*_{4n})$ (4) & N & N
     \\
     & $(D^*_{4m}/C_{2m}, D_{8n}^*/C_{4n})$ (11) & N & Y 
     \\
    \hline
    $(C_{2m}/C_{2m},T^*/T^*)$ (5)
     & $(C_{2m}/C_{2m},O^*/O^*)$ (7) & Y & N \\
     \small{$\gcd(m,6) = 1$}
     & $(C_{4m}/C_{4m}, T^*/T^*)$ (5) & N & N \\
     & $(D^*_{4m}/D^*_{4m},T^*/T^*)$ (14)  & N & N \\
     & $(C_{4m}/C_{2m}, O^*/T^*)$ (8) & N & N \\
     & $(D^*_{4m}/C_{2m}, O^*/T^*)$ (16) & N & Y \\
    \hline
    $(C_{6m}/C_{2m}, T^*/D_8^*)$ (6)
     &$(C_{12m}/C_{4m}, T^*/D_8^*)$ (6) &  N & N \\
     \small{$m$ odd and $3|n$}
     & $(D^*_{12m}/C_{2m}, O^*/D_8^*)$ (18) &  N & Y \\
    \hline
    $(C_{2m}/(C_{2m} \times O^*/O^*)$ (7)
     & $(C_{4m}/C_{4m},O^*/O^*)$ (7) & N & N\\
     \small{$\gcd(m,6) = 1$}
     & $(D_{4m}/D_{4m}, O^*/O^*)$ (15) &  N & Y \\
    \hline
    $(C_{2m}/C_{2m},I^*/I^*)$ (9)
     & $(C_{4m}/C_{4m},I^*/I^*)$ (9) &  N & N \\
      \small{$\gcd(m,30) = 1$}
     & $(D^*_{4m}/D^*_{4m},I^*/I^*)$ (19)  & N & Y \\
     \hline
    \hline
     \multicolumn{4}{|l|}{\small{$(*_1)$: if $\gcd(r,s)=1,$ $\gcd(m,n)=1$, $mn$ even or $r$ odd,}} \\
     \multicolumn{4}{|l|}{and $\gcd(n-sm,2mnr)= \gcd(n+sm,2mnr)$ } \\
     \hline
     \multicolumn{4}{|l|}{\small{$(*_2)$: if $\gcd(r,s)=1,$ $r$ even, $m$ odd, $n$ odd }}\\
     \multicolumn{4}{|l|}{\small{$\gcd(m,n)=1$ and $\gcd(n-sm,mnr)= \gcd(n+sm,mnr)$ }} \\
     \hline
\end{tabular}
\end{table}

\begin{table}
\caption{Geometric properties of involutions: small indices}
\label{small index quotient}

\begin{small}

\begin{tabular}{|l|l|c|c|}
	\hline
    group (family) & extension (family) & free & hyperelliptic\\
      &  & action &   \\
    \hline
    \hline
    $(C_{2m}/C_{2m},C_{2n}/C_{2n})$ (1)
     & $(C_{4m}/C_{4m},C_{2n}/C_{2n})$ (1) & 
     \small{$n$ odd}
     & N \\
   \small{$\gcd(n,m)=1$}
     & $(D^*_{4m}/D^*_{4m}, C_{2n}/ C_{2n})$ (2bis) & $n$ odd & N \\
     & $(C_{4m}/C_{2m}, C_{4n}/C_{2n})$ (1)  & $mn$ even & N \\
     & $(D^*_{4m}/C_{2m}, D^*_{4n}/C_{2n})$ (11) & N & Y \\
     & $(C_{4m}/C_{2m}, D^*_{4n}/C_{2n})$ (3)  &  $m$ even & N \\
     & $(C_{2m}/C_{2m},C_{4n}/C_{4n})$ (1)  &
 \small{$m$ odd}
 & N \\
     & $(C_{2m}/C_{2m},D^*_{4n}/D^*_{4n})$ (2) &  \small{$m$ odd} & N \\
     & $(D^*_{4m}/C_{2m}, C_{4n}/C_{2n})$ (3bis) & \small{$n$ even} & N \\
    \hline
    $(C_{4m}/C_{2m}, C_{4n}/C_{2n})$ (1)
     & $(C_{4m}/C_{4m},C_{4n}/C_{4n})$ (1) & N & N \\
    \small{$\gcd(m,n)=1$ and $mn$ even}
     & $(D^*_{8m}/D^*_{4m}, C_{4n}/C_{2n})$ (4bis)  & N & N \\
    
     & $(D_{8m}^*/C_{2m}, D_{8n}^*/C_{2n})$ (11)  & N & Y \\
     & $(C_{8m}/C_{2m}, C_{8n}/C_{2n})$ (1) & Y & N \\
     & $(C_{4m}/C_{2m}, D^*_{8n}/D^*_{4n})$ (4) & N & N \\
    \hline
    $(C_{2}/C_{2},C_{2n}/C_{2n})$ (1)
     & $(C_{4}/C_{4},C_{2n}/C_{2n})$ (1) &  \small{$n$ odd} & N \\
     & $(C_{2}/C_{2},C_{4n}/C_{4n})$ (1)  & Y & N \\
     & $(C_{2}/C_{2}, D^*_{4n}/D^*_{4n})$ (2)  & Y & N \\
     & $(C_{4}/C_{2}, C_{4n}/C_{2n})$ (1) &  \small{$n$ even} & N \\
     & $(C_{4}/C_{2}, D^*_{4n}/C_{2n})$ (3)  & N & Y \\
    \hline
    $(C_{2}/C_{2},C_{2}/C_{2})$ (1)
     & $(C_{4}/C_{2}, C_{4}/C_{2})$ (1) & N & Y \\
     & $(C_{4}/C_{4},C_{2}/C_{2})$ (1)  & Y & N \\
     & $(C_{2}/C_{2},C_{4}/C_{4})$ (1)  & Y & N \\
    \hline
    $(C_{2m}/C_{m}, C_{2n}/C_{n})$ ($1'$)
     & $(C_{2m}/C_{2m},C_{2n}/C_{2n})$ (1) & Y & N \\
    \small{$mn$ odd and $\gcd(m,n)=1$}
     & $(C_{4m}/C_{m}, C_{4n}/C_{n})$ ($1'$) & N & N \\
     & $(D^*_{4m}/C_{m}, D^*_{4n}/C_{n})$ ($11'$) & N & Y \\
     & $(C_{4m}/C_{m}, D^*_{4n}/C_{n})$ (34)  & N & N \\
     & $(D^*_{4m}/C_{m}, C_{4n}/C_{n})$ (34bis) & N & N \\
    \hline
    $(C_{2}/C_1, C_{2n}/C_{n})$ ($1'$)
     & $(C_{2}/C_{2},C_{2n}/C_{2n})$ (1) & Y & N \\
  \small{$n$ odd}
     & $(C_{4}/C_{1}, C_{4n}/C_{n})$ ($1'$)  & N & N \\
     & $(C_{4}/C_{1}, D^*_{4n}/C_{n})$ (34) & N & Y \\
    \hline
    $(C_{2}/C_{1}, C_{2}/C_{1})$ ($1'$)
     & $(C_{4}/C_{1}, C_{4}/C_{1})_1$ ($1'$) & N & Y \\
     & $(C_{2}/C_{2},C_{2}/C_{2})$ (1) & Y & N \\
    \hline
    $(C_{2}/C_{2},D_{4n}^*/D_{4n}^*)$ (2) 
     & $(C_{4}/C_{4}, D_{4n}^*/D_{4n}^*)$ (2)  & N & \small{$n$ even}\\
     & $(C_{2}/C_{2}, D_{8n^*}/D_{8n^*})$ (2)  & Y & N \\
     & $(C_{4}/C_{2}, D^*_{8n}/D^*_{4n})$ (4)  & N & \small{$n$ odd}\\
    \hline
    $(C_{2m}/C_{2m}, D_{8}/D_{8})$ (2) 
     & $(C_{4m}/C_{4m},D_{8}^*/D_{8}^*)$ (2)  & N & N \\
     \small{$m$ odd}
     & $(D^*_{4m}/D^*_{4m},D_{8}^*/D_{8}^*)$ (10)  & N & Y \\
     & $(C_{2m}/C_{2m},D_{16}^*/D_{16}^*)$ (2)  & Y & N \\
     & $(C_{4m}/C_{2m}, D^*_{16}/D^*_{8})$ (4)  & N & N \\
     & $(D^*_{4m}/C_{2m}, D^*_{16}/D^*_{8})$ (13bis) & N & N \\
    \hline
    $(C_{2}/C_{2},D_{8}/D_{8})$ (2)
     & $(C_{4}/C_{4}, D_{8}^*/D_{8}^*)$ (2)  & N & Y \\
     &  $(C_{2}/C_{2},D_{16}^*/D_{16}^*)$ (2)  & Y & N \\
     & $(C_{4}/C_{2}, D^*_{16}/D^*_{8})$ (4)  & N & N \\
    \hline
    $(C_2/C_2,T^*/T^*)$ (5)
     & $(C_4/C_4,T^*/T^*)$ (5) &  N & N \\
     & $(C_2/C_2,O^*,O^*)$ (7) &  Y & N \\
     & $(C_4/C_2, O^*/T^*)$ (8)  & N & Y \\
    \hline
    $(C_2/C_2, O^*/O^*)$ (7) & $(C_4/C_4,O^*/O^*)$ (7) &  N & Y \\
    \hline
    $(C_2/C_2,I^*/I^*)$ (9) & $(C_4/C_4,I^*/I^*)$ (9) & N & Y \\
    \hline
    
\end{tabular}
\end{small}

\end{table}

\begin{figure}[h!]\label{fig-montesinos}
\centering
\includegraphics[scale=0.7]{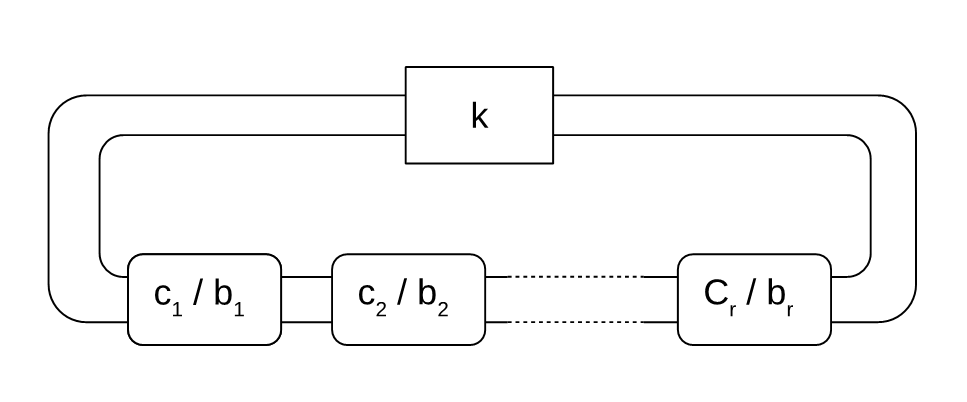}
\caption{Montesinos links}
\end{figure}

\subsection{More topological properties of the quotient orbifolds}\label{quotient}

Starting from the Seifert invariants of the quotient orbifold $\OO_f$, we can obtain other topological features of $\OO_f$.  

Seifert spherical orbifolds whose base orbifold is a sphere with at most two singular points are treated in \cite{MecchiaSeppi3}; in this case, the underlying topological space is a lens space that can be determined following \cite{geigeslange} or \cite{orlik}. If the base orbifold is a disk with at most one cone point an explicit description is given in \cite{dunbar2}. Also in this case the underlying topological space is a lens space. The underlying topological space of the other quotient orbifolds can be determined using the reduced local invariants, see Subsection~\ref{hyperelliptic}. It turns out that all the underlying topological spaces are spherical manifolds.

Finally, we focus on the singular set of the quotient orbifolds obtained by hyperelliptic involutions, so we suppose that the underlying topological space of $\OO_f$ is $S^3$. In many cases, the obtained Seifert fibration has the base orbifold $\BB_f$ of type $D^2(; b_1, b_2, \dots,b_r)$. This situation is completely described by \cite[Proposition 3]{dunbar}. Let $c_i/b_i$ be the local invariants of the corner points of $\BB_f$ and $e$ the Euler class of the fibration; the local invariants are defined mod integers so here we suppose that each $c_i/b_i$ is contained in the interval  $(-b_i/2,b_i/2)$. We define  $k$ in the following way:

$$k=-2e-\sum_i\frac{c_i}{b_i}$$ 

We recall that  $k$ is an integer. The singular set of  $\O_f$ is a Montesinos link as described in Figure~\ref{fig-montesinos}; each box labelled by a local invariant $c_i/b_i$ is a rational tangle defined by the rational number  $c_i/b_i$;  the box labelled by $k$ represents $k$ half twists.  More details can be found in \cite{dunbar}; we remark that we adopt here an opposite sign convention with respect to \cite{dunbar} for local invariants over cone points and corner points.

For the other quotient orbifolds, we obtain a fibration whose base orbifold is a sphere with some cone points. In this case, we can see the singular set as the union of some fibres of a fibration of the sphere. Also in this case, the singular set can be represented by a Montesinos link since each of these orbifolds admits another fibration whose base orbifold is of type $D^2(; b_1, b_2, \dots,b_r)$. The multiple fibrations of spherical orbifolds are described in \cite{mecchia-seppi2}. This completes the proof of Theorem~\ref{montesinos-link}.

\subsection*{Declarations}
The first author is a member of the national research group GNSAGA. The second author is a member of the Institut Camille Jordan, Université Claude Bernard Lyon 1. The authors have no conflicts of interest to declare relevant to this article's content. Our research did not generate or reuse research data.

\bibliographystyle{alpha}
\bibliography{references}

\begin{thebibliography}{McC02}

\bibitem[BS85]{bonahon-sibenmann}
F.~Bonahon and L.~Siebenmann.
\newblock The classification of {S}eifert fibred {$3$}-orbifolds.
\newblock In {\em Low-dimensional topology ({C}helwood {G}ate, 1982)},
  volume~95 of {\em London Math. Soc. Lecture Note Ser.}, pages 19--85.
  Cambridge Univ. Press, Cambridge, 1985.

\bibitem[BZ85]{burde-zieschang}
Gerhard Burde and Heiner Zieschang.
\newblock {\em Knots}, volume~5 of {\em De Gruyter Studies in Mathematics}.
\newblock Walter de Gruyter \& Co., Berlin, 1985.

\bibitem[CS03]{conway-smith}
John~H. Conway and Derek~A. Smith.
\newblock {\em On quaternions and octonions: their geometry, arithmetic, and
  symmetry}.
\newblock A K Peters, Ltd., Natick, MA, 2003.

\bibitem[Dun81]{dunbar2}
William~Dart Dunbar.
\newblock {\em F{ibered} {orbifolds} {and} {crystallographic} {groups}}.
\newblock ProQuest LLC, Ann Arbor, MI, 1981.
\newblock Thesis (Ph.D.)--Princeton University.

\bibitem[Dun88]{dunbar}
William~D. Dunbar.
\newblock Geometric orbifolds.
\newblock {\em Rev. Mat. Univ. Complut. Madrid}, 1(1-3):67--99, 1988.

\bibitem[DV64]{duval}
Patrick Du~Val.
\newblock {\em Homographies, quaternions and rotations}.
\newblock Oxford Mathematical Monographs. Clarendon Press, Oxford, 1964.

\bibitem[GL18]{geigeslange}
Hansj\"{o}rg Geiges and Christian Lange.
\newblock Seifert fibrations of lens spaces.
\newblock {\em Abh. Math. Semin. Univ. Hambg.}, 88(1):1--22, 2018.

\bibitem[HR85]{Hodgson-Rubinstein}
Craig Hodgson and J.~H. Rubinstein.
\newblock Involutions and isotopies of lens spaces.
\newblock In {\em Knot theory and manifolds ({V}ancouver, {B}.{C}., 1983)},
  volume 1144 of {\em Lecture Notes in Math.}, pages 60--96. Springer, Berlin,
  1985.

\bibitem[Kaw06]{Kawauchi}
Akio Kawauchi.
\newblock Topological imitations and {R}eni-{M}ecchia-{Z}immermann's
  conjecture.
\newblock {\em Kyungpook Math. J.}, 46(1):1--9, 2006.

\bibitem[Kim81]{kim}
Paik~Kee Kim.
\newblock Involutions on {K}lein spaces {$M(p,\,q)$}.
\newblock {\em Trans. Amer. Math. Soc.}, 268(2):377--409, 1981.

\bibitem[KM02]{kalliongis-miller}
John Kalliongis and Andy Miller.
\newblock Geometric group actions on lens spaces.
\newblock {\em Kyungpook Math. J.}, 42(2):313--344, 2002.

\bibitem[McC02]{mccullough}
Darryl McCullough.
\newblock Isometries of elliptic 3-manifolds.
\newblock {\em J. London Math. Soc. (2)}, 65(1):167--182, 2002.

\bibitem[Mon73]{montesinos}
Jos\'{e}~M. Montesinos.
\newblock Seifert manifolds that are ramified two-sheeted cyclic coverings.
\newblock {\em Bol. Soc. Mat. Mexicana (2)}, 18:1--32, 1973.

\bibitem[MS15]{mecchia-seppi}
Mattia Mecchia and Andrea Seppi.
\newblock Fibered spherical 3-orbifolds.
\newblock {\em Rev. Mat. Iberoam.}, 31(3):811--840, 2015.

\bibitem[MS19]{mecchia-seppi2}
Mattia Mecchia and Andrea Seppi.
\newblock Isometry groups and mapping class groups of spherical 3-orbifolds.
\newblock {\em Math. Z.}, 292(3-4):1291--1314, 2019.

\bibitem[MS20]{MecchiaSeppi3}
Mattia Mecchia and Andrea Seppi.
\newblock {On the diffeomorphism type of Seifert fibered spherical
  3-orbifolds}.
\newblock {\em {Rend. Ist. Mat. Univ. Trieste}}, 52:525--563, 2020.

\bibitem[Mye81]{myers}
Robert Myers.
\newblock Free involutions on lens spaces.
\newblock {\em Topology}, 20(3):313--318, 1981.

\bibitem[MZ04]{mecchia-zimmermann}
Mattia Mecchia and Bruno Zimmermann.
\newblock The number of knots and links with the same 2-fold branched covering.
\newblock {\em Q. J. Math.}, 55(1):69--76, 2004.

\bibitem[Orl72]{orlik}
Peter Orlik.
\newblock {\em Seifert manifolds}.
\newblock Lecture Notes in Mathematics, Vol. 291. Springer-Verlag, Berlin-New
  York, 1972.

\bibitem[Rub79]{Rubinstein}
J.~H. Rubinstein.
\newblock On {$3$}-manifolds that have finite fundamental group and contain
  {K}lein bottles.
\newblock {\em Trans. Amer. Math. Soc.}, 251:129--137, 1979.

\bibitem[Sei80]{seifert}
Herbert Seifert.
\newblock Topology of 3-dimensional fibered spaces.
\newblock In {\em A textbook of topology ({N}ew {Y}ork- {L}ondon, 1980)},
  volume~89 of {\em Pure Appl. Math.}, pages 139--152. Academic Press, New York
  - London, 1980.

\bibitem[TS31]{threlfall-seifert}
William Threlfall and Herbert Seifert.
\newblock Topologische {U}ntersuchung der {D}iskontinuit\"atsbereiche endlicher
  {B}ewegungsgruppen des dreidimensionalen sph\"arischen {R}aumes.
\newblock {\em Math. Ann.}, 104(1):1--70, 1931.

\bibitem[TS33]{threlfall-seifert2}
William Threlfall and Herbert Seifert.
\newblock Topologische {U}ntersuchung der {D}iskontinuit\"atsbereiche endlicher
  {B}ewegungsgruppen des dreidimensionalen sph\"arischen {R}aumes ({S}chlu\ss).
\newblock {\em Math. Ann.}, 107(1):543--586, 1933.

\bibitem[Vir76]{viro}
O.~Ja. Viro.
\newblock Nonprojecting isotopies and knots with homeomorphic coverings.
\newblock {\em Zap. Nau\v cn. Sem. Leningrad. Otdel. Mat. Inst. Steklov.
  (LOMI)}, 66:133--147, 207--208, 1976.
\newblock Studies in topology, II.

\end{thebibliography}

\end{document}